\documentclass[leqno]{article}
\usepackage{amsmath}
\usepackage{amsfonts,amssymb}
\usepackage[utf8]{inputenc}
\usepackage[T1]{fontenc}
\usepackage{hyperref}
\usepackage{amstext}
\usepackage{theorem}
\usepackage{enumerate}
\usepackage{float}
\usepackage{graphics,epsfig,times}
\usepackage{ae,aecompl}
\usepackage{mathrsfs}
\usepackage[all]{xy}
\usepackage[english]{babel}
\usepackage{fancyhdr}

\theoremstyle{definition}
\newtheorem{lemma}[equation]{Lemma}
\newtheorem{prop}[equation]{Proposition}
\newtheorem{cor}[equation]{Corollary}
\newtheorem{theo}[equation]{Theorem}
\newtheorem{remark}[equation]{Remark}

\newtheorem{definition}[equation]{Definition}
\numberwithin{equation}{section}
\newcommand{\OO}{\mathcal {O}}
\newcommand{\Q}{\mathcal{Q}}
\renewcommand{\P}{\mathcal{P}}
\newcommand{\N}{\mathcal{N}}

\newcommand{\R}{\mathcal{R}}
\newcommand{\T}{\mathcal{T}}
\newcommand{\bo}{\rm bo}
\newcommand{\Bo}{\rm Bo}
\newcommand{\p}{\mathfrak{p}}
\newcommand{\mfT}{\mathfrak T}
\newcommand{\dsp}{\displaystyle}
\newcommand{\RR}{\mathbb{R}}
\newcommand{\NN}{\mathbb{N}}
\newcommand{\mfQ}{\mathfrak{Q}}
\DeclareMathOperator*{\esssup}{ess\,sup}
\newcommand{\zetabar}{\underline{\zeta}}
\newcommand{\vbar}{\underline{v}}
\newcommand{\id}[1]{\left\vert_{\scriptstyle #1}\right.}
\newcommand{\nn}{\nonumber}

\newcommand{\demo}{\xqed{$\square$}}

\title{An improved result for the full justification of asymptotic models for the propagation of internal waves}

\author{Samer Israwi
\thanks{Laboratory of Mathematics-EDST and Faculty of Sciences I, Lebanese University, Beirut, Lebanon. ({\tt s\_israwi83@hotmail.com})({\tt rtalhouk@ul.edu.lb})}
\and Ralph Lteif \footnotemark[1]
 \thanks{LAMA, UMR 5127 CNRS, Université de Savoie Mont Blanc,
73376 Le Bourget du lac cedex, France. ({\tt ralphlteif\_90@hotmail.com}). SI, RL and RT are partially supported by the Lebanese University research
program (MAA group project).}
\and Raafat Talhouk \footnotemark[1]
} 

\begin{document}

\maketitle

\begin{abstract}
We consider here asymptotic models that describe the propagation of one-dimensional internal waves at the interface between two layers of immiscible fluids of different densities, under the rigid lid assumption and with uneven bottoms.
The aim of this paper is to show that the full justification result of the model obtained by Duch\^ene, Israwi and Talhouk [{\em SIAM J. Math. Anal.}, 47(1), 240–290], in the sense that it is consistent, well-posed, and that its solutions remain close to exact solutions of the full Euler system with corresponding initial data, can be improved in two directions.
The first direction is taking into account medium amplitude topography variations and the second direction is allowing strong nonlinearity using a new pseudo-symmetrizer, thus canceling out the smallness assumptions of the Camassa-Holm regime for the well-posedness and stability results.

\end{abstract}

\section{Introduction}
\subsection{Motivation}
The study of internal waves in a two-fluid system has attracted a lot of interests in a broad range of scientific communities because of the challenging
modeling, mathematical and numerical issues that arise in the analysis of this system. Such a configuration is commonly used in oceanography, where variations of temperature and salinity induce a density stratification.  Internal ocean waves, can interact with the bottom
topography and submerged structures as well as surface waves. Many authors have set a good theoretical background for
this problem in the flat bottom case. The case of uneven bottoms has been less investigated.
\\

The mathematical theory for this system of equation the so-called {\em full Euler system} derived in~\cite{Anh09},~\cite{BonaLannesSaut08} and~\cite{Duchene13} are fully nonlinear, and their direct study and computation remain a real
obstacle. In particular, the well-posedness of the equations in Sobolev space is challenging. Unlike the water wave problem (air-water interface), the Cauchy
problem associated with waves at the interface of two fluids of positive different densities over a flat bottom
is known to be ill-posed in Sobolev spaces in the absence of surface tension,
as Kelvin–Helmholtz instabilities appear. However, when adding a small amount of
surface tension, Lannes~\cite{Lannes13} proved that, thanks to a stability criterion, the problem
becomes well-posed with a time of existence that does not vanish as the surface tension
goes to zero and thus is consistent with the observations.
Consequently, a very large amount of work has been dedicated to the derivation of simplified,
asymptotic models, aiming at capturing the main characteristics of the 
flow with much simpler
equations, on the condition that the size of given parameters are small, thus reducing the framework to more
specific physical regimes, see Definition~1.1.
\\

Many models for a two-fluid system have already been derived and studied. Systems
under the rigid-lid assumption have first been investigated (see~\cite{Malcprimetseva89} or~\cite{Miyata85}, for
example).
Weakly and strongly nonlinear regimes have been derived
by Matsuno~\cite{Matsuno93}, Choi and Camassa~\cite{ChoiCamassa96,ChoiCamassa99}, generalizing the classical Green–Naghdi model (see~\cite{GreenNaghdi76}). Among other works, we would like to mention~\cite{BonaLannesSaut08} where many asymptotic models are presented and rigorously justified, in a wide range of regimes following a strategy initiated in~\cite{BonaChenSaut02,BonaChenSaut04}. In~\cite{GuyenneLannesSaut10}, the well-posedness and stability results have been proved for bi-fluidic shallow-water system without surface tension and under reasonable assumptions on the flow.
 Recently, Duch\^ene, Israwi and Talhouk derived in~\cite{DucheneIsrawiTalhouk13}, an asymptotic model for the propagation of one-dimensional internal
waves at the interface between two layers of immiscible 
fluids of different densities, under the
rigid lid assumption and with a 
flat bottom. They presented a new Green-Naghdi type model in
the Camassa-Holm (or medium amplitude) regime. They also proved that their model is fully justified, in the sense
that it is consistent, well-posed, and that its solutions remain close to exact solutions of the full
Euler system with corresponding initial data. Moreover, their system allows to fully justify any
well-posed and consistent lower order model; and in particular the so-called Constantin-Lannes
approximation (see~\cite{Duchene13}), which extends the classical Korteweg-de Vries equation in the Camassa-Holm
regime.\\ Let us note that all the aforementioned works are restricted to the case of a flat bottom, contrarily to the present work. As a matter of fact, topography plays an important role for internal waves. Indeed, in the ocean, internal waves are often propagating over variable topography, and this can cause quite dramatic deformation and transformation of internal waves. In this paper we extend the result of full justification obtained in~\cite{DucheneIsrawiTalhouk13} to the more complex case of non-flat topography, which is more reasonable in oceanography.
\\

Let us introduce some earlier results dealing with non-flat topography. A higher-order nonlinear
model for a two-layer
fluid of finite depth, is derived in~\cite{NachbinZarate09} to study the interaction of nonlinear internal waves with large
amplitude bottom topography that might vary rapidly over the characteristic
length scale of internal waves.
 Duch\^ene derived in~\cite{Duchene10} asymptotic models for the propagation of two and three dimensional gravity waves at the free surface and the interface between two layers of immiscible fluids
of different densities over an uneven bottom, following a method introduced by Bona, Lannes, and Saut in~\cite{BonaLannesSaut08} based on the expansion of the involved Dirichlet-to-Neumann operators. The same strategy was followed by Anh in~\cite{Anh09} to derive asymptotic models but taking into account surface tension, variable topography and for a large class of scaling regimes, furthermore the consistency of these asymptotic systems with the full Euler equations is established.
 To study the evolution of two-dimensional large amplitude internal waves in a two-layer system with
variable bottom topography, a new asymptotic model is derived by Barros and Choi in~\cite{ChoiBarros13}, the model is regularized to remove ill-posedness due to shear instability and it is further extended to include the effects of bottom topography and it is asymptotically equivalent to the strongly
nonlinear model proposed by Choi and Camassa~\cite{ChoiCamassa99}. Finally, a bi-fluidic Green-Naghdi model  that allows a non-flat topography, and horizontal dimension $d = 2$ is derived and justified (in the sense of consistency) in~\cite{DucheneIsrawiTalhouk14}. Its derivation follows directly from classical results concerning the one-layer case. However, the previously mentioned works are restricted to the formal level. In fact, the models derived in these papers are systematically justified
by a consistency result, and do not provide the full justification, contrarily to the present work.

\subsection{Outline of the paper}
This paper deals with the propagation of one-dimensional internal waves between two layers of fluids of different densities over an uneven bottom  under the following assumptions: the fluids are homogeneous, immiscible, ideal, incompressible, irrotationnal, and under the only influence of gravity.
We assume that the surface is confined by a flat rigid lid. The domain of the two layers is infinite in the horizontal space variable.

The main idea of our paper is to improve the
result of full justification of the Green-Naghdi type model in the Camassa-Holm (or medium amplitude) regime obtained by Duch\^ene, Israwi and Talhouk in~\cite{DucheneIsrawiTalhouk13} in two directions. Since the aforementioned work is restricted to the case of a flat topography, the first direction is to relax this assumption by taking 
into account medium amplitude topography variations, which is more reasonable in oceanography. To this end, we introduce a new parameter $\beta$ to characterize the shape of the bottom. Moreover we assume that $\beta=\OO(\sqrt \mu)$ which corresponds to the physical case of a bottom with medium variations in amplitude. In the linear analysis of the model obtained in~\cite{DucheneIsrawiTalhouk13}, the only place where the smallness assumption of the Camassa-Holm regime, $\epsilon=\OO(\sqrt\mu)$ is used apart from the simplifications it offers when constructing the model, stands in the proof of the energy estimates~\cite[Lemma 6.5]{DucheneIsrawiTalhouk13}. As a matter of fact, the second direction of our improvement is to construct a new pseudo-symmetrizer that allows to cancel the use of the Camassa-Holm regime assumptions ($\epsilon=\OO(\sqrt \mu)$ and $\beta=\OO(\sqrt \mu)$) from the energy estimate proof.  Therefore these assumptions can be relaxed for the well-posedness and stability results, regardless of their necessity for the derivation of the model. 
 Our new model possesses a structure similar to symmetrizable quasilinear systems that allows the study of the properties, and in particular energy
estimates, for the linearized system, thus
allowing its full justification, following the classical theory of hyperbolic systems. We prove that the full Euler system is consistent with our model, and
that our system is well-posed (in the sense of Hadamard) in Sobolev spaces, and stable with respect
to perturbations of the equations.
Finally, we prove that the solutions of our system approach the solutions of the full Euler system, with as good a precision as $\mu$  is small.
 \subsection{Organization of the paper}
  The present paper is organized as follows. In Section 2  we introduce the non-dimensionalized full Euler system and the Green-Naghdi model. Section 3 is dedicated to the construction of the new asymptotic model. Sections 4 and 5 contain the necessary ingredients for the proof of our results. In Section 6, we explain the full justification of asymptotic models and we state its main ingredients.
 \medskip
 
 \paragraph{Notations}

 In the following, $C_0$ denotes any nonnegative constant whose exact expression is of no importance.\\
 The notation $a\lesssim b$ means that
 $a\leq C_0\ b$ and we write $A=\OO(B)$ if $A\leq C_0\ B$.\\
 We denote by $C(\lambda_1, \lambda_2,\dots)$ a nonnegative constant depending on the parameters
 $\lambda_1$, $\lambda_2$,\dots and whose dependence on the $\lambda_j$ is always assumed to be nondecreasing.\\
 We use the condensed notation
\[
A_s=B_s+\left\langle C_s\right\rangle_{s>\underline{s}},
\]
to express that $A_s=B_s$ if $s\leq \underline{s}$ and $A_s=B_s+C_s$ if $s> \underline{s}$.\\
 Let $p$ be any constant
 with $1\leq p< \infty$ and denote $L^p=L^p(\RR)$ the space of all Lebesgue-measurable functions
 $f$ with the standard norm \[\vert f \vert_{L^p}=\big(\int_{\RR}\vert f(x)\vert^p dx\big)^{1/p}<\infty.\] The real inner product of any functions $f_1$
 and $f_2$ in the Hilbert space $L^2(\RR)$ is denoted by
\[
 (f_1,f_2)=\int_{\RR}f_1(x)f_2(x) dx.
 \]
 The space $L^\infty=L^\infty(\RR)$ consists of all essentially bounded, Lebesgue-measurable functions
 $f$ with the norm
\[
 \vert f\vert_{L^\infty}= \esssup \vert f(x)\vert<\infty.
\]
Let $k \in \NN$, we denote by $W^{k,\infty}=W^{k,\infty}(\RR)=\{f \in L^{\infty}, |f|_{W^{k,\infty}}< \infty\}$, where $ |f|_{W^{k,\infty}}=\displaystyle{\sum_{\alpha \in \NN , \alpha\leq k} |\partial_x^{\alpha} f|_{L^{\infty}}}$.
 For any real constant $s\geq0$, $H^s=H^s(\RR)$ denotes the Sobolev space of all tempered
 distributions $f$ with the norm $\vert f\vert_{H^s}=\vert \Lambda^s f\vert_{L^2} < \infty$, where $\Lambda$
 is the pseudo-differential operator $\Lambda=(1-\partial_x^2)^{1/2}$.\\
For a given $\mu>0$, we denote by $H^{s+1}_\mu(\RR)$ the space $H^{s+1}(\RR)$ endowed with the norm
\[\big\vert\ \cdot\ \big\vert_{H^{s+1}_\mu}^2 \ \equiv \ \big\vert\ \cdot\ \big\vert_{H^{s}}^2\ + \ \mu \big\vert\ \cdot\ \big\vert_{H^{s+1}}^2\ .\]
 For any function $u=u(t,x)$ and $v(t,x)$
 defined on $ [0,T)\times \RR$ with $T>0$, we denote the inner product, the $L^p$-norm and especially
 the $L^2$-norm, as well as the Sobolev norm,
 with respect to the spatial variable $x$, by $(u,v)=(u(t,\cdot),v(t,\cdot))$, $\vert u \vert_{L^p}=\vert u(t,\cdot)\vert_{L^p}$,
 $\vert u \vert_{L^2}=\vert u(t,\cdot)\vert_{L^2}$ , and $ \vert u \vert_{H^s}=\vert u(t,\cdot)\vert_{H^s}$, respectively.\\
 We denote $L^\infty([0,T);H^s(\RR))$ the space of functions such that $u(t,\cdot)$ is controlled in $H^s$, uniformly for $t\in[0,T)$:
 \[\big\Vert u\big\Vert_{L^\infty([0,T);H^s(\RR))} \ = \ \esssup_{t\in[0,T)} \vert u(t,\cdot)\vert_{H^s} \ < \ \infty.\]
 Finally, $C^k(\RR)$ denote the space of
 $k$-times continuously differentiable functions.\\
 For any closed operator $T$ defined on a Banach space $X$ of functions, the commutator $[T,f]$ is defined
  by $[T,f]g=T(fg)-fT(g)$ with $f$, $g$ and $fg$ belonging to the domain of $T$. The same notation is used for $f$ an operator mapping the domain of $T$ into itself.\\
  
We conclude this section by the nomenclature that we use to describe the different regimes
that appear in the present work. A regime is defined through restrictions on the set of admissible
dimensionless parameters of the system, which are precisely defined in~\eqref{regimes}, below.
\begin{definition}{\em (Regimes)}\label{defregime}.
We define the so-called {\em shallow water regime} for two layers of comparable depths the set of parameters:
\begin{multline} \label{eqn:defRegimeSWmr}
\P_{\rm SW} \ \equiv \ \Big\{ (\mu,\epsilon,\delta,\gamma,\beta,\bo):\ 0\ < \ \mu \ \leq \ \mu_{\max}, \ 0 \ \leq \ \epsilon \ \leq \ 1, \ \delta \in (\delta_{\min},\delta_{\max}),  \Big.\\
 \Big. \ 0\ \leq \ \gamma\ <\ 1, \ 0\leq\beta\leq\beta_{\max},\  \bo_{\min}\leq \bo\leq \infty \ \Big\},
\end{multline}
with given $0\leq \mu_{\max},\delta^{-1}_{\min},\delta_{\max},\bo_{\min}^{-1},\beta_{\max}<\infty$.\\
\\
The additional key restrictions for the validity of the model \eqref{eq:Serre2} are as follows:
\begin{equation} \label{eqn:defRegimeCHmr}
\P_{\rm CH}  \equiv  \left\{ (\mu,\epsilon,\delta,\gamma,\beta,\bo) \in \P_{\rm SW}:\ \epsilon  \leq  M \sqrt{\mu} ,\quad \beta\leq  M \sqrt{\mu}  \quad \mbox{ and } \quad \nu \equiv  \frac{1+\gamma\delta}{3\delta(\gamma+\delta)}-\frac{1}{\bo} \ge  \nu_{0}  \ \right\},
\end{equation}
with given $0\leq M,\nu_0^{-1}<\infty$.\\
\\
We denote for convenience
\[ M_{\rm SW} \ \equiv \ \max\big\{\mu_{\max},\delta_{\min}^{-1},\delta_{\max},\bo_{\min}^{-1},\beta_{\max}\big\}, \quad M_{\rm CH} \ \equiv \ \max\big\{M_{\rm SW},M,\nu_0^{-1}\big\}.\]
 \end{definition}

\section{Previously obtained models}

\subsection{The full Euler system}
We briefly recall the governing
equations of a two-layer flow in the aforementioned configuration, that we call {\em full Euler system} (using non-dimensionalized variables and the Zakharov/Craig-Sulem formulation)~\cite{CraigSulem93,Zakharov68}. We let the
reader refer to~\cite{Anh09,BonaLannesSaut08,Duchene13} for more details.

\begin{equation}\label{eqn:EulerCompletAdim}
\left\{ \begin{array}{l}
\displaystyle\partial_{ t}{\zeta} \ -\ \frac{1}{\mu}G^{\mu}\psi\ =\ 0,  \\ \\
\displaystyle\partial_{ t}\Big(H^{\mu,\delta}\psi-\gamma \partial_x{\psi} \Big)\ + \ (\gamma+\delta)\partial_x{\zeta} \ + \ \frac{\epsilon}{2} \partial_x\Big(|H^{\mu,\delta}\psi|^2 -\gamma |\partial_x {\psi}|^2 \Big) \\ \hspace{5cm} = \mu\epsilon\partial_x\N^{\mu,\delta}-\mu\frac{\gamma+\delta}{\bo}\frac{\partial_x \big(k(\epsilon\sqrt\mu\zeta)\big)}{{\epsilon\sqrt\mu}} \ ,
\end{array}
\right.
\end{equation}
where we denote
\[  \N^{\mu,\delta} \ \equiv \ \frac{\big(\frac{1}{\mu}G^{\mu}\psi+\epsilon(\partial_x{\zeta})H^{\mu,\delta}\psi \big)^2\ -\ \gamma\big(\frac{1}{\mu}G^{\mu}\psi+\epsilon(\partial_x{\zeta})(\partial_x{\psi}) \big)^2}{2(1+\mu|\epsilon\partial_x{\zeta}|^2)}.
      \]
$\zeta(t,x)$ represent the deformation of the interface between the two layers and $b(x)$ represent the deformation of the bottom,  $\psi$ is the trace of the velocity potential of the upper-fluid at the interface.\\
\\
The function $ k(\zeta)=-\partial_x \Big(\frac1{\sqrt{1+|\partial_x\zeta|^2}}\partial_x\zeta\Big)$ denotes the mean curvature of the interface.

Let $a$(resp. $a_b$) be the maximal vertical deformation of the interface(resp. bottom) with respect to its rest position. We denote by $\lambda$ a characteristic horizontal length, say the wavelength of the interface; $d_1$ (resp. $d_2$) the depth of the upper (resp. lower) layer; and $\rho_1$ (resp. $\rho_2$) is the density of the upper (resp. lower) layer, $g$ the gravitational acceleration, $\sigma$ the interfacial tension coefficient. In the following we use $\bo=\mu\Bo$ instead of {\em the classical Bond number}, $\Bo$ which measures the ratio of gravity forces over capillary forces. Consequently we introduce the following dimensionless parameters
\begin{equation}\label{regimes}
\mu \ = \ \frac{d_1^2}{\lambda^2}, \quad \epsilon \ = \ \frac{a}{d_1}, \quad \beta \ = \ \frac{a_b}{d_1},\quad \delta \ = \ \frac{d_1}{d_2},\quad \gamma \ = \ \frac{\rho_1}{\rho_2},\quad \bo\ =\ \frac{g(\rho_2-\rho_1)d_{1}^2}{\sigma},\end{equation}
where $\epsilon$ (resp. $\beta$) measures the amplitude of the deformation at the interface (resp. bottom) with respect to the depth of the upper layer of fluid, and $\mu$ is the shallowness parameter.\\

Finally, $G^{\mu}$ and $H^{\mu,\delta}$ are the so-called Dirichlet-Neumann operators, defined as follows:
\begin{definition}[Dirichlet-Neumann operators]
Let $\zeta\in H^{t_0+1}(\RR)$, $t_0>1/2$, such that there exists $h_0>0$ with
$h_1 \ \equiv\  1-\epsilon\zeta \geq h_0>0$ and $h_2 \ \equiv \ \frac1\delta +\epsilon \zeta-\beta b\geq h_0>0$, and let $\psi\in L^2_{\rm loc}(\RR),\ \partial_x \psi\in H^{1/2}(\RR)$.
 Then we define
\begin{eqnarray*}
G^{\mu}\psi  &\equiv & G^{\mu}[\epsilon\zeta]\psi  \equiv  \sqrt{1+\mu|\epsilon\partial_x\zeta|^2}\big(\partial_n \phi_1 \big)\id{z=\epsilon\zeta}  =  -\mu\epsilon(\partial_x\zeta) (\partial_x\phi_1)\id{z=\epsilon\zeta}+(\partial_z\phi_1)\id{z=\epsilon\zeta},\\
H^{\mu,\delta}\psi  &\equiv  & H^{\mu,\delta}[\epsilon\zeta,\beta b]\psi \equiv  \partial_x \big(\phi_2\id{z=\epsilon\zeta}\big)  =  (\partial_x\phi_2)\id{z=\epsilon\zeta}+\epsilon(\partial_x \zeta)(\partial_z\phi_2)\id{z=\epsilon\zeta},
\end{eqnarray*}
where $\phi_1$ and $\phi_2$ are uniquely defined (up to a constant for $\phi_2$) as the solutions in $H^2(\RR)$ of the Laplace's problems:
\begin{eqnarray}
\label{eqn:Laplace1} &&\left\{
\begin{array}{ll}\left(\ \mu\partial_x^2 \ +\  \partial_z^2\ \right)\ \phi_1=0 & \mbox{ in } \Omega_1\equiv \{(x,z)\in \RR^{2},\ \epsilon{\zeta}(x)<z<1\}, \\
\partial_z \phi_1 =0  & \mbox{ on } \Gamma_{\rm t}\equiv \{(x,z)\in \RR^{2},\ z=1\}, \\
\phi_1 =\psi & \mbox{ on } \Gamma\equiv \{(x,z)\in \RR^{2},\ z=\epsilon \zeta\},
\end{array}
\right.\\ 
\label{eqn:Laplace2}&&\left\{
\begin{array}{ll}
\left(\ \mu\partial_x^2\ + \ \partial_z^2\ \right)\ \phi_2=0 & \mbox{ in } \Omega_2\equiv\{(x,z)\in \RR^{2},\ -\frac{1}{\delta}+\beta b(x)<z<\epsilon\zeta\}, \\
\partial_{n}\phi_2 = \partial_{n}\phi_1 & \mbox{ on } \Gamma,\\
\partial_{n}\phi_2 =0 & \mbox{ on } \Gamma_{\rm b}\equiv \{(x,z)\in \RR^{2},\ z=-\frac{1}{\delta}+\beta b(x)\}.
\end{array}
\right.
\end{eqnarray}
\end{definition}

The existence and uniqueness of a solution to \eqref{eqn:Laplace1}-\eqref{eqn:Laplace2}, and therefore the well-posedness of the Dirichlet-Neumann operators follow from classical arguments detailed, for example, in~\cite{Lannes}.

\subsection{The Green-Naghdi model}
In the following, we construct Green-Naghdi type model for the system~\eqref{eqn:EulerCompletAdim}, that is asymptotic
model with precision $\OO(\mu^2)$, in the sense of consistency.
The key ingredient for constructing shallow water asymptotic models comes from the expansion of the Dirichlet-Neumann operators given in~\cite{Duchene10,DucheneIsrawiTalhouk14}, with respect to the shallowness parameter, $\mu$. When replacing the Dirichlet-to-Neumann operators by their truncated expansion, and after
straightforward computations, one is able to deduce the Green-Naghdi model, that we disclose
below. This model has been introduced in~\cite{Duchene13} (with a flat bottom) and generalized in~\cite{DucheneIsrawiTalhouk14}. It is also convenient to introduce a new velocity variable, namely the shear mean velocity $v$ is equivalently defined as \begin{equation}\label{defv}v\equiv  u_2 - \gamma  u_1\end{equation} where $u_1$ and $u_2$ are the horizontal velocities integrated across each layer:\\
$ u_1(t,x)  =  \frac{1}{h_1(t,x)}\int_{\epsilon\zeta(t,x)}^{1} \partial_x \phi_1(t,x,z) \ dz$ and $ u_2(t,x)  =  \frac{1}{h_2(t,x)}\int_{-\frac1\delta+\beta b(x)}^{\epsilon\zeta(t,x)} \partial_x \phi_2(t,x,z) \ dz$, where $\phi_1$ and $\phi_2$ are the solutions to the Laplace's problems \eqref{eqn:Laplace1}-\eqref{eqn:Laplace2}.\\
\\
The expansions of the Dirichlet-Neumann operators may be written in terms of the new variable $v$.\\
\\
Plugging the expansions given in~\cite{DucheneIsrawiTalhouk14} into the full Euler system \eqref{eqn:EulerCompletAdim}, and withdrawing all $\OO(\mu^2)$ terms yields  ( in the unidimensional case $d=1$),
\begin{equation}\label{eq1}
\left\{ \begin{array}{l}
\displaystyle\partial_{ t}{\zeta} \ + \ \partial_x \Big(\frac{h_1h_2}{h_1+\gamma h_2} v\Big)\ =\ 0,  \\ \\
\displaystyle\partial_{ t}\Bigg(  v \ + \ \mu\overline{\Q}[h_1,h_2] v \Bigg) \ + \ (\gamma+\delta)\partial_x{\zeta} \ + \ \frac{\epsilon}{2} \partial_x\Big(\dfrac{h_1^2 -\gamma h_2^2 }{(h_1+\gamma h_2)^2}  v^2\Big) \ = \ \qquad \\
\hfill \displaystyle\mu\epsilon\partial_x\big(\overline{\R}[h_1,h_2] v  \big) +\mu\frac{\gamma+\delta}{\bo}\partial_{x}^3\zeta,
\end{array}
\right.
\end{equation}
where we denote $h_1=1-\epsilon\zeta$ and $h_2=\delta^{-1}+\epsilon\zeta-\beta b$, as well as

\begin{eqnarray*}\overline{\Q}[h_1,h_2]v&=&\T[h_2,\beta b]\Big(\frac{h_1v}{h_1+\gamma h_2}\Big)-\gamma\T[h_1,0]\Big(\frac{-h_2v}{h_1+\gamma h_2}\Big),\\
&=&-\dfrac{1}{3h_2}\partial_x\Big(h_2^3\partial_x\big(\dfrac{h_1v}{h_1+\gamma h_2}\big)\Big)+\dfrac{1}{2h_2}\beta\Big[\partial_x\Big(h_2^2(\partial_x b) \dfrac{h_1 v}{h_1+\gamma h_2}\Big)-h_2^2(\partial_x b)\partial_x \big(\dfrac{h_1v}{h_1+\gamma h_2}\big)\Big]\\& \qquad + & \beta^2(\partial_x b)^2\big(\dfrac{h_1v}{h_1+\gamma h_2}\big)-\gamma\Big[\dfrac{1}{3h_1}\partial_x\Big(h_1^3\partial_x\big(\dfrac{h_2v}{h_1+\gamma h_2}\big)\Big)\Big].\end{eqnarray*}

\begin{eqnarray*}\overline{\R}[h_1,h_2] v&=& \dfrac{1}{2}\Big(-h_2\partial_x(\frac{h_1v}{h_1+\gamma h_2})+\beta(\partial_x{b})(\frac{h_1v}{h_1+\gamma h_2})\Big)^2-\dfrac{\gamma}{2}\Big(h_1\partial_x(\frac{-h_2v}{h_1+\gamma h_2})\Big)^2\\& \qquad - & (\frac{h_1v}{h_1+\gamma h_2})\T[h_2,\beta b]\Big(\frac{h_1v}{h_1+\gamma h_2}\Big)+\gamma(\frac{-h_2v}{h_1+\gamma h_2})\T[h_1,0]\Big(\frac{-h_2v}{h_1+\gamma h_2}\Big),\end{eqnarray*}
with  \begin{equation*}\T[h,b]V\equiv \dfrac{-1}{3h}\partial_x(h^3\partial_x V)+ \dfrac{1}{2h}[\partial_x(h^2(\partial_x b) V)-h^2(\partial_x b)(\partial_x V)]+(\partial_x b)^2 V .\end{equation*}
\section{Construction of the new model}
In the following section, we will use several additional transformations, in order to produce an equivalent model (in the sense of consistency) which possesses a structure similar to symmetrizable quasilinear systems that allows the study of the subsequent sections.\\ 
The present work is limited to the so-called Camassa-Holm regime, that is
using two additional assumptions $\epsilon=\OO(\sqrt\mu)$ and  $\beta=\OO(\sqrt\mu)$ (deformation of the interface and the one of the bottom are of medium amplitude). We manipulate the Green-Naghdi system~\eqref{eq1}, systematically withdrawing $\OO(\mu^2,\mu\epsilon^2,\mu\beta^2,\mu\epsilon\beta)$ terms, in order to recover our model.\\
One can check that the following approximations formally hold:
 \begin{eqnarray*}
        \overline{\Q}[h_1,h_{2}] v  &=&  - {\lambda}\partial_x^2 v-\epsilon\frac{\gamma+\delta}3 \left((\theta-\alpha)v \partial_x^2\zeta +  (\alpha+2\theta)\partial_x(\zeta\partial_x v)-\theta\zeta\partial_x^2 v\right)\\&\qquad+&\beta\frac{\gamma+\delta}3 \left((\frac{\alpha_1}{2}+\theta_1)v \partial_x^2 b +  (\alpha_1+2\theta_1)\partial_x(b\partial_x v)-\theta_1b\partial_x^2 v\right)  \\&\qquad+&  \OO(\epsilon^2,\beta^2,\epsilon\beta), \\
        \overline{\R}[h_1,h_{2}] v   &=&  \alpha\left(\frac12 (\partial_x  v)^2+\frac13 v\partial_x^2 v\right) + \OO(\epsilon,\beta).
 \end{eqnarray*}
with \begin{equation}\label{eq3}
{\lambda}=\frac{1+\gamma\delta}{3\delta(\gamma+\delta)}\ , \quad \alpha=\dfrac{1-\gamma}{(\gamma+\delta)^2} \quad \mbox{ and } \quad \theta=\dfrac{(1+\gamma\delta)(\delta^2-\gamma)}{\delta(\gamma+\delta)^3}\ ,
\end{equation}
and \begin{equation}\label{eq3}
\alpha_1=\dfrac{1}{(\gamma+\delta)^2} \quad \mbox{ and } \quad \theta_1=\dfrac{\delta(1+\gamma\delta)}{(\gamma+\delta)^3}\ .
\end{equation}

Plugging these expansions into system~\eqref{eq1} and using the same techniques as in~\cite[Section 4.2]{DucheneIsrawiTalhouk13} but with a different symmetric operator ${\mathfrak T}[\epsilon\zeta,\beta b]$ defined below, yields a simplified model, with the same order of precision of the original model (that is $\OO(\mu^2)$) in the Camassa-Holm regime.\\
Let us first introduce the operator ${\mathfrak T}[\epsilon\zeta,\beta b]$.
 \begin{equation}\label{defT}
      {\mathfrak T}[\epsilon\zeta,\beta b]V \ = \ q_1(\epsilon\zeta,\beta b)V \ - \ \mu\nu \partial_x\Big(q_2(\epsilon\zeta,\beta b)\partial_xV \Big),
      \end{equation}
      with $q_i(X,Y)\equiv 1+\kappa_i X +\omega_i Y $ ($i=1,2$) and $\nu,\kappa_1,\kappa_2,\omega_1,\omega_2,\varsigma$ are constants to be determined.
\\
\\
The first order  $\OO(\mu)$ terms may be canceled with a proper choice of $\nu$, making use of the fact that the second equation of the Green-Naghdi model~\eqref{eq1} yields
\begin{equation*}
\partial_t v=-(\gamma+\delta) \partial_x \zeta -\dfrac{\epsilon}{2} \partial_x \Big(\dfrac{\delta^2-\gamma}{(\gamma+\delta)^2}|v^2|\Big)+\OO(\epsilon^2,\mu,\epsilon\beta).\end{equation*}
Indeed it follows that
\begin{equation}\label{eq123}
\dfrac{\gamma+\delta}{\bo}\partial_x^3 \zeta=-\dfrac{1}{\bo} \partial_x^2\partial_t v -\dfrac{3\epsilon}{2\bo}\dfrac{\delta^2-\gamma}{(\gamma+\delta)^2}\partial_x ((\partial_x v)^2)-\dfrac{\epsilon}{\bo}\dfrac{\delta^2-\gamma}{(\gamma+\delta)^2} v \partial_x^3 v+\OO(\epsilon^2,\mu,\epsilon\beta).\end{equation}
Using again that~\eqref{eq1} yields $\partial_t v=-(\gamma+\delta) \partial_x \zeta + \OO(\epsilon, \beta, \mu)$ and 
$\partial_t \zeta=\frac{-1}{\gamma+\delta}\partial_x v + \OO(\epsilon, \beta, \mu)$, it becomes clear, that one can adjust $\kappa_1$, $\kappa_2$, $\omega_1$, $\omega_2 \in \RR$ so that all terms involving $\zeta$ and its derivatives are withdrawn.\\
\\
In order to deal with ($v\partial_x^3 v$) terms, $\varsigma \in \RR$ is to be determined. In fact, these terms  appear after replacing  the term $\dfrac{\gamma+\delta}{\bo}\partial_x^3 \zeta$ of the second equation of~\eqref{eq1} by its expression given in~\eqref{eq123}. Thus one defines the constants $\nu,\kappa_1,\kappa_2,\omega_1,\omega_2,\varsigma$ as follow:

         \begin{equation}\label{deftnu}
 \nu \ = \ \lambda-\frac1{\bo} \ = \ \frac{1+\gamma\delta}{3\delta(\gamma+\delta)}-\frac1{\bo},
      \end{equation}
\begin{equation}\label{defkappa}
 (\lambda-\frac1{\bo})\kappa_1 \ = \ \frac{\gamma+\delta}{3}(2\theta-\alpha) , \quad   (\lambda-\frac1{\bo})\kappa_2 \ = \ (\gamma+\delta)\theta , 
  \end{equation}
  \begin{equation}\label{defomega}
 (\lambda-\frac1{\bo})\omega_1 \ = \ -\theta_1\frac{
  (\gamma+\delta)}{3} , \quad   (\lambda-\frac1{\bo})\omega_2 \ = \ -\dfrac{(\gamma+\delta)}{3}(\alpha_1+2\theta_1) ,
  \end{equation}

      \begin{equation}\label{defvarsigma}
  (\lambda-\frac1{\bo})\varsigma \ = \ \frac{2\alpha-\theta}{3} \ - \ \frac{1}{\bo} \dfrac{\delta^2 -\gamma }{(\delta+\gamma )^2}.
      \end{equation}

Finally, one can check that
\begin{multline}\label{eqcons}
       {\mathfrak T}[\epsilon\zeta,\beta b](\partial_t  v +\epsilon \varsigma  v\partial_x  v) - q_1(\epsilon\zeta,\beta b) \partial_{ t}\Big(  v \ + \ \mu\overline{\Q}[h_1,h_{2}] v\Big) +\mu q_1(\epsilon\zeta,\beta b)  \Big(\frac{\gamma+\delta}{\bo}\partial_x^3\zeta+\epsilon\partial_{ x} \big(\overline{\R}[h_1,h_{2}]  v\big)  \Big) \\
         =  \epsilon \varsigma q_1(\epsilon\zeta,\beta b)  v\partial_x  v   -  \mu\epsilon \frac{2\alpha}3\partial_x\big((\partial_x   v)^2\big)+\mu\beta \omega (\partial_x \zeta)(\partial_x^2 b) +\OO(\mu^2,\mu\epsilon^2,\mu\beta^2,\mu\epsilon\beta)\end{multline}

            where we denote $\omega=\dfrac{(\gamma+\delta)^2}{3}\big(\dfrac{\alpha_1}{2}+\theta_1\big)$.\\
            \\
       When plugging the estimate~\eqref{eqcons} in~\eqref{eq1}, and after multiplying the second equation by $q_1(\epsilon\zeta,\beta b)$, we obtain the following system of equations:
      \begin{equation}\label{eq:Serre2}\left\{ \begin{array}{l}
      \partial_{ t}\zeta +\partial_x\left(\dfrac{h_1 h_{2}}{h_1+\gamma h_2}  v\right)\ =\  0,\\ \\
      \mathfrak T[\epsilon\zeta,\beta b] \left( \partial_{ t}   v + \epsilon\varsigma { v } \partial_x {  v } \right) + (\gamma+\delta)q_1(\epsilon\zeta,\beta b)\partial_x
      \zeta +\frac\epsilon2 q_1(\epsilon\zeta,\beta b) \partial_x \left(\frac{h_1^2  -\gamma h_2^2 }{(h_1+\gamma h_2)^2}| v|^2-\varsigma | v|^2\right)\\ \qquad \qquad \qquad \ = \    -  \mu \epsilon\frac23\alpha \partial_x\big((\partial_x  v)^2\big)+\mu \beta \omega(\partial_x \zeta)(\partial_x^2 b) ,
      \end{array} \right. \end{equation}

      \medskip

\begin{prop}[Consistency]\label{th:ConsSerreVar}
For $\p= (\mu,\epsilon,\delta,\gamma,\beta,\bo) \in \P_{\rm SW}$, let $U^\p=(\zeta^\p,\psi^\p)^\top$ be a family of solutions of the full Euler system~\eqref{eqn:EulerCompletAdim}  such that there exists $T>0$, $s\geq s_0+1$, $s_0 > 1/2$   for which $(\zeta^\p,\partial_x\psi^\p)^\top$ is bounded in $L^{\infty}([0,T); H^{s+N})^2$ with  sufficiently large $N$, and uniformly with respect to $\p \in \P_{\rm SW}$.
 Moreover assume that $b\in H^{s+N}$ and there exists $h_{01}>0$ such that
\begin{equation*}
h_1=1-\epsilon\zeta^\p\geq h_{01}>0,  \quad h_2=1/\delta +\epsilon\zeta^\p-\beta b\geq h_{01}>0.
\end{equation*}

Define $v^\p$ as in \eqref{defv}. Then $(\zeta^\p, v^\p)^\top$ satisfies \eqref{eq:Serre2} up to a remainder term, $R=(0,r)^\top$, bounded by $$\|r\|_{L^{\infty}([0,T); H^{s})}\leq(\mu^2+\mu\epsilon^2+\mu\beta^2+\mu\epsilon\beta) C,$$
with $C=C(h_{01}^{-1},M_{\rm SW},|b|_{H^{s+N}},\|(\zeta^\p,\partial_x{\psi^\p})^\top\|_{L^{\infty}([0,T); H^{s+N})^2}).$
\end{prop}
\emph{Proof}.\\
 Let $U=(\zeta,\psi)^\top$ satisfy the hypothesis above withdrawing the explicit dependence with respect to parameters $\p$ for the sake of readability. We know from~\cite[Proposition 3.14]{DucheneIsrawiTalhouk14} that $(\zeta, v)^\top$ satisfies the system \eqref{eq1} up to a remainder $R_0=(0,r_0)^\top$ bounded by,
 $$\|r_0\|_{L^{\infty}([0,T); H^{s})}\leq\mu^2 C_1,$$
with $C_1=C(h_{01}^{-1},M_{\rm SW},|b|_{H^{s+N}},\|(\zeta^\p,\partial_x{\psi^\p})^\top\|_{L^{\infty}([0,T); H^{s+N})^2})$, uniformly with respect to $\p \in\P_{\rm SW}$.\\
 The proof now consists in checking that all terms neglected in the above calculations can be rigorously estimated in the same way. We do not develop each estimate, but rather provide the precise bound on
the various remainder terms. One has
  \begin{multline*}
 \Big\vert    \partial_t\big(   \overline{\Q}[h_1,h_{2}]v \big)    -   \big[- \lambda\partial_x^2\partial_t  v-\epsilon\frac{\gamma+\delta}3 \partial_t\left((\beta-\alpha) v \partial_x^2\zeta  +  (\alpha+2\beta)\partial_x(\zeta\partial_x v)-\beta\zeta\partial_x^2 v\right)\\+\beta\frac{\gamma+\delta}3 \partial_t \left((\frac{\alpha_1}{2}+\theta_1)v \partial_x^2 b +  (\alpha_1+2\theta_1)\partial_x(b\partial_x v)-\theta_1b\partial_x^2 v\right) \big] \Big\vert_{H^{s}}\\
  \leq (\epsilon^2+\beta^2+\epsilon \beta) C(s+3), \ \end{multline*}
with $C(s+3)\equiv C\Big(M_{\rm SW},h_{01}^{-1},\big\vert \zeta \big\vert_{H^{s+3}},\big\vert \partial_t\zeta  \big\vert_{H^{s+2}},\big\vert  v \big\vert_{H^{s+3}},\big\vert  \partial_t v \big\vert_{H^{s+2}},\big\vert  b \big\vert_{H^{s+3}}\Big)$, and
\[
    \Big\vert     \partial_x\big(  \overline{\R}[h_1,h_{2}] v\big)   -  \partial_x\big[\alpha\big(\frac12 (\partial_x  v)^2+\frac13 v\partial_x^2  v\big)\big]\Big\vert_{H^{s}} \\
    \leq
   ( \epsilon+\beta) C(s+3) .\
\]
Then, since $(\zeta, v)$ satisfies~\eqref{eq1}, up to the remainder $R_0$, one has
\[\big\vert  \partial_t  v+(\gamma+\delta)\partial_x\zeta  \big\vert_{H^{s}}+  \big\vert\partial_t\zeta+\frac{1}{\gamma+\delta}\partial_x  v  \big\vert_{H^{s}}\ \leq \ (\epsilon+\beta) C(s+3) + \big\vert R_0 \big\vert_{H^{s}} .\]  
It follows that \eqref{eqcons} is valid up to a remainder $R_1$, bounded by
$$|R_1|_{H^s}\leq (\mu^2+\mu\epsilon^2+\mu\beta^2+\mu\epsilon\beta)C(s+3)+\mu(\epsilon+\beta+\mu)|R_0|_{H^s}$$
Finally, $(\zeta, v)$ satisfies \eqref{eq:Serre2} up to a remainder $R$, bounded by
$$|R|_{H^s}\leq|R_0+R_1|_{H^s}\leq(\mu^2+\mu\epsilon^2+\mu\beta^2+\mu\epsilon\beta)C.$$
where we use that
\[\big\vert  v  \big\vert_{H^{s+3}}+\big\vert \partial_t  v \big\vert_{H^{s+2}}\leq C.\]
The estimate on $v$ follows directly from the identity $\partial_x \Big(\frac{h_1h_2}{h_1+\gamma h_2}v\Big)  =  -\frac1\mu G^{\mu,\epsilon}\psi=\partial_t\zeta$. The estimate on $\partial_t v$ can be proved, for example, following~\cite[Prop. 2.12]{Duchene10}.
This concludes the proof of Proposition ~\ref{th:ConsSerreVar}.

\demo
\section{Preliminary results}

In this section, we recall the operator $\mfT[\epsilon\zeta,\beta b]$, defined in~\eqref{defT}:
\begin{equation}
{\mathfrak T}[\epsilon\zeta,\beta b]V \ = \ q_1(\epsilon\zeta,\beta b)V \ - \ \mu\nu \partial_x\Big(q_2(\epsilon\zeta,\beta b)\partial_xV \Big).
\end{equation}
with $\nu,\kappa_1,\kappa_2,\omega_1,\omega_2$ are constants and $\nu=\dfrac{1+\gamma\delta}{3\delta(\gamma+\delta)}-\dfrac{1}{\bo}\geq\nu_0>0.$
\\
The operator  $\mfT[\epsilon\zeta,\beta b]$, has exactly the same structure as the one introduced in~\cite{DucheneIsrawiTalhouk13} but it depends also on the deformation of the bottom and plays an important role in the energy estimate and the local well-posednees of the system~\eqref{eq:Serre2}.\\
In the following, we seek sufficient conditions to ensure the strong ellipticity of the operator  ${\mfT}$ which will yield to the well-posedness and continuity of the inverse ${\mfT}^{-1}$.
As a matter of fact, this condition, namely~\eqref{CondEllipticity} (and similarly the classical non-zero depth condition,~\eqref{CondDepth}) simply consists in assuming that the deformation of the interface is not too large as given in~\cite{DucheneIsrawiTalhouk13} but here we have to take into account the topographic variation that plays a role in~\eqref{CondDepth} and in~\eqref{CondEllipticity}.
For fixed $\zeta\in L^\infty$ and $b\in L^\infty$, the restriction reduces to an estimate on $\epsilon_{\max}\big\vert \zeta\big\vert_{L^\infty}+\beta_{\max}|b|_{L^\infty}$ with $\epsilon_{\max},\beta_{\max}=\min(M\sqrt{\mu_{\max}},1)$, and~\eqref{CondDepth}-\eqref{CondEllipticity} hold uniformly with respect to $(\mu,\epsilon,\delta,\gamma,\beta,\bo)\in\P_{\rm CH}.$
\\
\\
Let us recall the non-zero depth condition
\begin{equation}\label{CondDepth}\tag{H1}
\exists \  h_{01}>0, \quad \mbox{such that} \quad \inf_{x\in \RR} h_1\geq h_{01}>0,\quad
\inf_{x\in \RR} h_2\geq h_{01}>0.
\end{equation}
where $h_1=1-\epsilon\zeta$ and $h_2=1/\delta+\epsilon\zeta-\beta b$ are the depth of respectively the upper and the lower layer of the fluid. \\
\\
It is straightforward to check that, since for all $(\mu,\epsilon,\delta,\gamma,\beta,\bo)\in\P_{\rm CH}$, the following condition
\[\epsilon_{max}|\zeta|_{L^\infty}+\beta_{max}|b|_{L^\infty}< \min(1,\dfrac{1}{\delta_{max}})\]
is sufficient to define $h_{01}>0$ such that \eqref{CondDepth} is valid independently of $(\mu,\epsilon,\delta,\gamma,\beta,\bo)\in\P_{\rm CH}.$\\
In the same way, we introduce the condition
\begin{equation}\label{CondEllipticity}\tag{H2}
\exists\ h_{02}>0 , \mbox{ such that } \quad  \inf_{x\in \RR} \left(1+\epsilon\kappa_1\zeta+\beta\omega_1 b \right) \ge  h_{02} \ > \ 0 \ ; \qquad    \inf_{x\in \RR}  \left( 1+\epsilon \kappa_2\zeta +\beta\omega_2 b \right)\geq h_{02} \ > \ 0.
\end{equation}

{\em In what follows, we will always assume that~\eqref{CondDepth} and~\eqref{CondEllipticity} are satisfied}. It is a consequence of our work that such assumption may be imposed only on the initial data, and then is automatically satisfied over the relevant time scale.\\
\\
Now the preliminary results proved in~\cite[Section 5]{DucheneIsrawiTalhouk13} remain valid for the operator ${\mathfrak T}[\epsilon\zeta,\beta b]$.
Before asserting the strong ellipticity of the operator  ${\mfT}$, let us
first recall the quantity $\vert \cdot \vert_{H^{1}_\mu}$, which is defined as
\[\forall \  v\in H^1(\RR),\ \quad
\vert\ v\ \vert^2_{H^{1}_\mu}\ =\ \vert\ v\ \vert^2_{L^2}\ +\ \mu\ \vert\ \partial_x v\ \vert^2_{L^2},
\]
 and is equivalent to the $H^1(\RR)$-norm but not uniformly with respect to $\mu$.
 We define by, $H^1_\mu(\RR)$ the space $H^1(\RR)$ endowed with this norm and $(H^1_\mu(\RR))^\star$ the space $H^{-1}(\RR)$ the dual space of  $H_\mu^1(\RR)$.\\
 
The following lemma gives an important invertibility result on $\mfT$.

\begin{lemma}\label{Lem:mfT}
Let $(\mu,\epsilon,\delta,\gamma,\beta,\bo)\in\P_{\rm CH}$ and $\zeta \in L^{\infty}(\RR)$, $b\in L^{\infty}(\RR)$ such that~\eqref{CondEllipticity} is satisfied.
 Then the operator
\[
{\mfT}[\epsilon\zeta,\beta b]: H^1_\mu(\RR)\longrightarrow (H^1_\mu(\RR))^\star
\]
is uniformly continuous and coercive. More precisely, there exists $c_0>0$ such that
\begin{eqnarray}
({\mfT} u,v)  \ & \leq& \  c_0\vert u\vert_{H^1_\mu}\vert v\vert_{H^1_\mu}  ; \label{continous} \\
({\mfT} u,u) \ & \geq & \frac1{c_0}\vert u\vert_{H^1_\mu}^{2} \label{coercive}
\end{eqnarray}
with $c_0=C(M_{\rm CH},h_{02}^{-1},\epsilon\big\vert \zeta\big\vert_{L^\infty},\beta\big\vert b \big\vert_{L^\infty} )$.
\medskip

Moreover, the following estimates hold:
Let $ s_0>\frac{1}{2}$ and $s\geq0$,

\begin{enumerate}
\item[(i)]  If $\zeta, b \in H^{s_0}(\RR)\cap H^{s} (\RR)$ and $u\in  H^{s+1}(\RR)$ and $v\in H^{1}(\RR)$, then:
\begin{equation}
\label{eq:estLambdaT}\big\vert \big(\Lambda^s \mfT[\epsilon\zeta,\beta b]  u, v \big) \big\vert \ \leq \  C_0 \left((\epsilon|\zeta|_{H^{s_0}}+\beta|b|_{H^{s_0}})\big\vert u\big\vert_{H^{s+1}_\mu} +\big\langle (\epsilon|\zeta|_{H^{s}}+\beta|b|_{H^{s}}) \big\vert u\big\vert_{H^{s_0+1}_\mu} \big\rangle_{s>s_0}\right) \big\vert v\big\vert_{H^{1}_\mu} . \ 
\end{equation}

\item[(ii)]  If $\zeta, b \in H^{s_0+1}\cap H^{s} (\RR)$ , $u\in H^{s}(\RR)$ and $v\in H^{1}(\RR)$, then:
\begin{align}
\label{eq:estComT}\big\vert \big( \big[\Lambda^s, \mfT[\epsilon\zeta,\beta b]\big]u, v\big) \big\vert &\leq 
C_0 \Big((\epsilon|\zeta|_{H^{s_0+1}}+\beta|b|_{H^{s_0+1}})\big\vert  u\big\vert_{H^{s}_\mu}\nn\\ & \qquad \qquad \qquad + \big\langle (\epsilon|\zeta|_{H^{s}}+\beta|b|_{H^{s}})\big\vert u\big\vert_{H^{s_0+1}_\mu}\big\rangle_{s>s_0+1}\Big) \big\vert v\big\vert_{H^{1}_\mu}\nn\\ &\leq \max(\epsilon,\beta)
C_0 \Big((|\zeta|_{H^{s_0+1}}+|b|_{H^{s_0+1}})\big\vert  u\big\vert_{H^{s}_\mu}\nn\\ & \qquad \qquad \qquad + \big\langle (|\zeta|_{H^{s}}+|b|_{H^{s}})\big\vert u\big\vert_{H^{s_0+1}_\mu}\big\rangle_{s>s_0+1}\Big) \big\vert v\big\vert_{H^{1}_\mu},\nn\\
\end{align}
\end{enumerate}
where $C_0=C(M_{\rm CH},h_{02}^{-1})$.
\end{lemma}
The following lemma then gives some properties of the inverse operator ${\mfT}^{-1}$.
\begin{lemma}\label{proprim}
Let $(\mu,\epsilon,\delta,\gamma,\beta,\bo)\in\P_{\rm CH}$ and $\zeta \in L^{\infty}(\RR)$, $b\in L^{\infty}(\RR)$ such that~\eqref{CondEllipticity} is satisfied.
 Then the operator
\[
{\mfT}[\epsilon\zeta,\beta b]: H^2(\RR)\longrightarrow L^2(\RR)
\]
is one-to-one and onto. Moreover, one has the following estimates:
\begin{enumerate}
\item[(i)] $({\mfT}[\epsilon\zeta,\beta b])^{-1}:L^2\to H^1_\mu(\RR)$ is continuous. More precisely, one has
\[
\parallel {\mfT}^{-1}\parallel_{L^2(\RR)\rightarrow H^1_\mu(\RR)} \ \leq \  c_0,
\]
with $c_0=C(M_{\rm CH},h_{02}^{-1}, \epsilon\big\vert \zeta\big\vert_{L^\infty} ,\beta \big\vert b \big\vert_{L^\infty} )$.
\item[(ii)] Additionally, if $\zeta, b \in H^{s_0+1}(\RR)$ with $s_0>\frac{1}{2}$, then one has for any $0< s\leq s_0 +1$,
\[
\parallel {\mfT}^{-1}\parallel_{H^s(\RR)\rightarrow H^{s+1}_\mu(\RR)}\ \leq \ c_{s_0+1}.
\]
\item[(iii)] If $\zeta, b \in H^{s}(\RR)$ with $ s\ge s_0 +1,\; s_0>\frac{1}{2}$, then one has
\[
\parallel {\mfT}^{-1}\parallel_{H^s(\RR)\rightarrow H^{s+1}_\mu(\RR)}\ \leq \ c_{s}
\]
\end{enumerate}
where $c_{\bar s}=C(M_{\rm CH},h_{02}^{-1},\epsilon|\zeta|_{H^{\bar s}},\beta|b|_{H^{\bar s}})$, thus uniform with respect to $(\mu,\epsilon,\delta,\gamma,\beta,\bo)\in\P_{\rm CH}$.
\end{lemma}

Finally, let us introduce the following technical estimate, which is used several times in the subsequent sections.
\begin{cor}\label{col:comwithT}
Let $(\mu,\epsilon,\delta,\gamma,\beta,\bo)\in\P_{\rm CH}$ and $\zeta, b \in H^{s}(\RR)$ with $ s\ge s_0 +1,\; s_0>\frac{1}{2}$, such that~\eqref{CondEllipticity} is satisfied. Assume moreover that $u\in H^{s-1}(\RR)$ and that $v\in H^{1}(\RR)$. Then one has
\begin{eqnarray}
\big|\big( \ \big[\Lambda^s, \mfT^{-1}[\epsilon\zeta,\beta b]\big] u\ ,\ \mfT[\epsilon\zeta,\beta b]  v\ \big)\big|  \ &= & \ \big|\big(\  \big[\Lambda^s,\mfT[\epsilon\zeta,\beta b] \big] {\mfT}^{-1}  u\ ,\  v\ \big)\big| \nn \\
& \leq&  \ \max(\epsilon,\beta)\ C(M_{\rm CH},h_{02}^{-1},\big\vert \zeta \big\vert_{H^s},\big\vert b\big\vert_{H^s})\big\vert u \big\vert_{H^{s-1}} \big\vert v \big\vert_{H^{1}_\mu}
\label{eq:comwithT}\end{eqnarray}
\end{cor}

\section{Linear analysis}
This section is devoted to the study of the properties and in particular energy estimates, for the linearized system associated to our nonlinear asymptotic system~\eqref{eq:Serre2}, following the classical theory of quasilinear hyperbolic systems. More precisely, we establish new linear estimates independent of the Camassa-Holm assumptions $\epsilon=\OO(\sqrt \mu)$ and $\beta=\OO(\sqrt \mu)$. To establish these estimates we propose a new pseudo-symmetrizer of the system that allows us to cancel the use of the smallness assumptions of the Camassa-Holm regime in the proof of Lemma~5.18, in fact the assumption on the deformation of the interface $\epsilon=\OO(\sqrt \mu)$ was necessary for the proof of~\cite[Lemma 6.5]{DucheneIsrawiTalhouk13}. Therefore, these assumptions can be relaxed for the well-posedness and stability results (Theorem~6.1 and Theorem~6.6 respectively), regardless of their necessity for the derivation of the model~\eqref{eq:Serre2}.
\\
\\
Let us recall the system~\eqref{eq:Serre2}.
 \begin{equation}\label{eq:Serre2var}\left\{ \begin{array}{l}
      \partial_{ t}\zeta +\partial_x\left(\dfrac{h_1 h_{2}}{h_1+\gamma h_2}  v\right)\ =\  0,\\ \\
      \mathfrak T[\epsilon\zeta,\beta b] \left( \partial_{ t}   v + \epsilon\varsigma { v } \partial_x {  v } \right) + (\gamma+\delta)q_1(\epsilon\zeta,\beta b)\partial_x
      \zeta +\frac\epsilon2 q_1(\epsilon\zeta,\beta b) \partial_x \left(\frac{h_1^2  -\gamma h_2^2 }{(h_1+\gamma h_2)^2}| v|^2-\varsigma | v|^2\right)\\ \qquad \qquad \qquad \ = \    -  \mu \epsilon\frac23\alpha \partial_x\big((\partial_x  v)^2\big)+\mu \beta \omega(\partial_x \zeta)(\partial_x^2 b) ,
      \end{array} \right. \end{equation}

with $h_1=1-\epsilon \zeta$ , $h_2=1/\delta+\epsilon \zeta-\beta b$ , $q_i(X,Y)=1+\kappa_i X+\omega_i Y$ ($i=1,2$) , $\kappa_i,\omega_i,\varsigma$ defined in \eqref{defkappa},\eqref{defomega},\eqref{defvarsigma}, and
\[
\mfT[\epsilon\zeta,\beta b] V= q_1(\epsilon\zeta,\beta b)V -\dsp \mu \nu\partial_x \left(q_2(\epsilon\zeta,\beta b)\partial_x V \right).
\]
In order to ease the reading, we define the function
\[ f:X\to\frac{(1-X)(\delta^{-1}+X-\beta b)}{1-X+\gamma(\delta^{-1}+X-\beta b)},\] and \[ g:X\to\Big(\frac{(1-X)}{1-X+\gamma(\delta^{-1}+X-\beta b)}\Big)^{2}.\]
One can easily check that
\[
f(\epsilon\zeta) \ =\ \dsp\frac{h_1h_2}{h_1+\gamma h_2}, \qquad
f'(\epsilon\zeta) \ =\ \dsp\frac{h_1^2-\gamma h_2^2}{(h_1+\gamma h_2)^2}\quad\mbox{ and }\quad g(\epsilon\zeta)=\Big(\frac{h_1}{h_1+\gamma h_2}\Big)^2 .\]
Additionally, let us denote \[
\kappa=\frac23 \alpha=\frac23\frac{1-\gamma}{(\delta+\gamma)^2}\quad \mbox{ and}
\quad q_3(\epsilon\zeta)=\frac12\big(\frac{h_1^2-\gamma h_2^2}{(h_1+\gamma h_2)^2}-\varsigma\big),\]
so that one can rewrite,
\begin{equation}\label{eqn:Serre2varf}\left\{ \begin{array}{l}
\dsp \partial_{ t}\zeta +f(\epsilon\zeta)\partial_x  v+\epsilon \partial_x\zeta f'(\epsilon\zeta)  v -\beta\partial_x b g(\epsilon\zeta)v \ =\  0,\\ \\
\dsp {\mathfrak T} \left( \partial_{ t}  v + \frac{\epsilon}{2} \varsigma \partial_x({v }^2) \right) + (\gamma+\delta)q_1(\epsilon\zeta,\beta b)\partial_x \zeta + \epsilon q_1(\epsilon\zeta,\beta b)\partial_x(q_3(\epsilon\zeta)  {v}^2) \\ \hspace{7cm} =- \mu\epsilon\kappa\partial_x\big((\partial_x v)^2\big)+\mu \beta \omega (\partial_x \zeta) (\partial_x^2 b).
\end{array} \right. \end{equation}
with $\partial_x(q_3(\epsilon\zeta))=\dfrac{-\gamma\epsilon\partial_x \zeta (h_1+h_2)^2 +\gamma \beta \partial_x b h_1(h_1+h_2)}{(h_1+\gamma h_2)^3}.$
\\
\\
The equations can be written after applying ${\mathfrak T}^{-1}$ to the second equation
in~\eqref{eqn:Serre2varf} as
\begin{equation}\label{condensedeq}
\partial_tU+A[U]\partial_xU +B[U]\ = \ 0,\end{equation}
 with
\begin{equation}\label{defA0A1}
A[U]=\begin{pmatrix}
\epsilon f'(\epsilon \zeta) v&f(\epsilon\zeta)\\
\mfT^{-1}(Q_0(\epsilon\zeta,\beta b )\cdot+\epsilon^2 Q_1(\epsilon\zeta,\beta b,v) \cdot)& \epsilon \mfT^{-1}(\mfQ[\epsilon\zeta,\beta b,v] \cdot)+\epsilon\varsigma v
\end{pmatrix}
\end{equation}
\begin{equation}\label{defB}
B[U]=\begin{pmatrix}
 -\beta\partial_x b g(\epsilon\zeta)v\\
\epsilon\mfT^{-1}\Big(\dfrac{\gamma\beta q_1(\epsilon\zeta,\beta b)h_1(h_1+h_2)v^2}{(h_1+\gamma h_2)^3}\partial_x b\Big)
\end{pmatrix}
\end{equation}
where $Q_0(\epsilon\zeta,\beta b),Q_1(\epsilon\zeta,\beta b,v)$ are defined as

\begin{equation}\label{defQ0Q1}
Q_0(\epsilon\zeta,\beta b) \ = \ (\gamma+\delta)q_1(\epsilon\zeta,\beta b)-\mu\beta\omega \partial_x^2 b,\quad Q_1(\epsilon\zeta,\beta b,v)=-\gamma q_1(\epsilon \zeta,\beta b)\dfrac{(h_1+h_2)^2}{(h_1+\gamma h_2)^3}{ v}^2
\end{equation}
and the operator $\mfQ[\epsilon,\beta b,v]$ defined by
\begin{equation}\label{defmfQ}
\mfQ[\epsilon\zeta,\beta b,v]f \ \equiv \ 2q_1(\epsilon\zeta,\beta b)q_3(\epsilon\zeta)v f +\mu\kappa \partial_x(f \partial_x v).
\end{equation}
The following sections are devoted to the proof of energy estimates for the following
initial value problem around some reference state
$\underline{U}=(\underline{\zeta},\underline{v})^\top$:
\begin{equation}\label{SSlsys}
	\left\lbrace
	\begin{array}{l}
	\dsp\partial_t U+A[\underline{U}]\partial_x U+B[\underline{U}]=0;
        \\
	\dsp U_{\vert_{t=0}}=U_0.
	\end{array}\right.
\end{equation}

\subsection{Energy space}

Let us first remark that by construction, the new pseudo-symmetrizer is given by
\begin{equation}\label{defZ}Z[\underline{U}]=\begin{pmatrix}
 \frac{Q_0(\epsilon\underline{\zeta},\beta b)+\epsilon^2 Q_1(\epsilon\underline{\zeta},\beta b,\underline{v})}{f(\epsilon\underline{\zeta})}& 0 \\
0&\mfT[\epsilon\underline{\zeta},\beta b]
\end{pmatrix}
\end{equation}
Note however that one should add an additional assumption in order to ensure that our pseudo-symmetrizer is defined and positive which is:
\begin{equation}\label{H3} \tag{H3} \exists \  h_{03} >0 \  \  \mbox {such that} \  \    Q_0(\epsilon\underline{\zeta},\beta b)+\epsilon^2 Q_1(\epsilon\underline{\zeta},\beta b,\underline{v})\geq h_{03} > 0.\end{equation}

We define now the $X^s$ spaces, which are the energy spaces for this problem.
\begin{definition}\label{defispace}
 For given $s\ge 0$ and $\mu,T>0$, we denote by $X^s$ the vector space $H^s(\RR)\times H^{s+1}_\mu(\RR)$ endowed with the norm
\[
\forall\; U=(\zeta,v) \in X^s, \quad \vert U\vert^2_{X^s}\equiv \vert \zeta\vert^2 _{H^s}+\vert v\vert^2 _{H^s}+ \mu\vert \partial_xv\vert^2 _{H^s},
\]
while $X^s_T$ stands for the space of $U=(\zeta,v)$ such that $U\in C^0([0,\frac{T}{\max(\epsilon,\beta)}];X^{s})$, and $\partial_t U \in  L^\infty([0,\frac{T}{\max(\epsilon,\beta)}]\times \RR)$, endowed with the canonical norm
\[
 \Vert U\Vert_{X^s_T}\equiv \sup_{t\in [0,T/\max(\epsilon,\beta)]}\vert U(t,\cdot)\vert_{X^s}+\esssup_{t\in [0,T/\max(\epsilon,\beta)],x\in\RR}\vert \partial_t U (t,x)\vert.
\]
\end{definition}

A natural energy for the initial value problem \eqref{SSlsys} is now given by:
\begin{equation}
 E^s(U)^2=(\Lambda^sU,Z[\underline{U}]\Lambda^sU)=( \Lambda^s\zeta,\frac{Q_0(\epsilon\underline{\zeta},\beta b)+\epsilon^2 Q_1(\epsilon\underline{\zeta},\beta b,\underline{v})}{f(\epsilon\underline{\zeta})}\Lambda^s\zeta)+\left(\Lambda^s v,\mfT[\epsilon\underline{\zeta},\beta b] \Lambda^s v\right).
\end{equation}
In order to ensure the equivalency of $X^s$ with the energy of the pseudo-symmetrizer it requires to add the additional assumption given in \eqref{H3}.\\
\\
The link between $E^s(U)$ and the $X^s$-norm is investigated in the following Lemma. 
\begin{lemma}\label{lemmaes}
 Let $\p=(\mu, \epsilon, \delta, \gamma, \beta, \bo)\in \P_{\rm CH}$, $s\geq 0$ , $ \underline{U}\in L^{\infty}(\RR)$ and $b \in W^{2,\infty}(\RR)$, satisfying~\eqref{CondDepth}, \eqref{CondEllipticity}, and~\eqref{H3}. Then
$E^s(U)$ is equivalent to the $\vert \cdot\vert_{X^s}$-norm.\\ More precisely, there exists $c_0=C(M_{\rm CH},h_{01}^{-1},h_{02}^{-1},h_{03}^{-1}, \epsilon|\underline{U}|_{L^\infty},\beta|b|_{W^{2,\infty}})>0$ such that
\[
\frac1{c_0}E^s(U) \ \leq \ \big\vert U \big\vert_{X^s} \ \leq \  c_0 E^s(U).
\]
\end{lemma}
\emph{Proof}.\\
This is a straightforward application of Lemma~\ref{Lem:mfT}, and that for $Q_0(\epsilon\underline{\zeta},\beta b)+\epsilon^2 Q_1(\epsilon\underline{\zeta},\beta b,\underline{v})\geq h_{03}>0$ and $f(\epsilon\underline{\zeta})>0$,
\[
\inf_{x\in\RR} \frac{Q_0(\epsilon\underline{\zeta},\beta b)+\epsilon^2 Q_1(\epsilon\underline{\zeta},\beta b,\underline{v})}{f(\epsilon\underline{\zeta})}\geq  \inf_{x\in\RR} \Big( Q_0(\epsilon\underline{\zeta},\beta b))+\epsilon^2 Q_1(\epsilon\underline{\zeta},\beta b,\underline{v})\Big)\Big(\sup_{x\in\RR} f(\epsilon\underline{\zeta})\Big)^{-1},\]
\[
\sup_{x\in\RR} \frac{Q_0(\epsilon\underline{\zeta},\beta b)+\epsilon^2 Q_1(\epsilon\underline{\zeta},\beta b,\underline{v})}{f(\epsilon\underline{\zeta})}\leq \sup_{x\in\RR} \Big(Q_0(\epsilon\underline{\zeta},\beta b)+\epsilon^2 Q_1(\epsilon\underline{\zeta},\beta b,\underline{v})\Big)\Big(\inf_{x\in\RR} f(\epsilon\underline{\zeta})\Big)^{-1}.\]
where we recall that if \eqref{CondDepth} is satisfied then, $\underline{h_1}=1-\epsilon\zetabar, \underline{h_2}=1/\delta+\epsilon\zetabar-\beta b$ satisfy
\[\inf_{x\in\RR}\underline{h_1}\geq h_{01}, \quad \sup_{x\in\RR}|\underline{h_1}|\leq 1 + 1/\delta,  \quad \inf_{x\in\RR}\underline{h_2}\geq h_{01}, \quad \sup_{x\in\RR}|\underline{h_2}|\leq1+1/\delta.\]\demo

\begin{lemma}\label{lem-nsZ}
  Let $\p=(\mu, \epsilon, \delta, \gamma, \beta, \bo)\in \P_{\rm CH}$, and let $U=(\zeta_u,u)^\top \in L^\infty$,  $b \in W^{2,\infty}$ satisfies~\eqref{CondDepth},\eqref{CondEllipticity} and \eqref{H3}. Then for any
 $V,W\in X^{0}$, one has
\begin{equation}\label{eq:est-Z}
\Big\vert \ \Big(\ Z[U] V \ , \ W \ \Big) \ \Big\vert \  \leq \ \ C \ \big\vert V\big\vert_{X^0}\big\vert W\big\vert_{X^0},
\end{equation}
with $C=C(M_{\rm CH},h_{01}^{-1},h_{02}^{-1},\epsilon\big\vert U\big\vert_{L^\infty},\beta \big\vert b \big\vert_{W^{2,\infty}})\ $.
\\
\\
Moreover, if $U \in X^s, b \in H^{s+2} ,V\in X^{s-1}$ with $s\geq s_0+1,\ s_0>1/2$, then one has
\begin{eqnarray}
\label{eq:est-com-Z}
\Big\vert  \Big(\ \big[ \Lambda^s,Z[U]\big] V \ , \ W \ \Big) \Big\vert  \ &\leq& \ C\ \big\vert V\big\vert_{X^{s-1}}\big\vert W\big\vert_{X^0}\\
\label{eq:est-com-Z-Z}
\Big\vert  \Big(\ \big[ \Lambda^s,Z^{-1}[U]\big] V \ , \ Z[U] W \ \Big) \Big\vert  \ &\leq& \  C\ \big\vert V\big\vert_{H^{s-1}\times H^{s-1}}\big\vert W\big\vert_{X^0}
\end{eqnarray}
with $C=C(M_{\rm CH},h_{01}^{-1},h_{02}^{-1},\epsilon\big\vert U \big\vert_{X^s},\beta  \big\vert b \big\vert_{H^{s+2}})\ $.
\end{lemma}
\emph{Proof}.\\
The Lemma~5.13 is proved using Cauchy-Schwarz inequality, Lemma~4.2 and Corollary~4.8. We do not detail the proof, and refer to~\cite[Lemma 6.4]{DucheneIsrawiTalhouk13}.\demo

\subsection{Energy estimates}

Our aim is to establish {\em a priori} energy estimates concerning our linear system. In order to be able to use the linear analysis to both the well-posedness and stability of the nonlinear system, we consider the following modified system
\begin{equation}\label{SSlsysm}
	\left\lbrace
	\begin{array}{l}
	\dsp\partial_t U+A[\underline{U}]\partial_x U+B[\underline{U}]= F;
        \\
	\dsp U_{\vert_{t=0}}=U_0.
	\end{array}\right.
\end{equation}
where we added a right-hand-side $F$, whose properties will be precised in the following Lemmas.\\
\\
We begin by asserting a basic $X^0$ energy estimate, that we extend to $X^s$ space ($s>3/2$) later on.

\begin{lemma}[$X^0$ energy estimate]\label{Lem:L2z}
Set $(\mu, \epsilon, \delta, \gamma, \beta, \bo)\in \P_{\rm CH}$.\\ Let $T>0$, $s_0>1/2$ and $U\in L^\infty ( [0,T/\max(\epsilon,\beta)];X^0)$ and $\underline{U},\partial_x \underline{U}\in L^\infty([0,T/\max(\epsilon,\beta)]\times \RR)$ and $b \in H^{s_0+3}$ such that $\partial_t \underline{U} \in L^\infty([0,T/\max(\epsilon,\beta)]\times \RR)$ and $\underline{U}, b $ satisfies~\eqref{CondDepth},\eqref{CondEllipticity}, and \eqref{H3} and $U,\underline{U}$ satisfy system~\eqref{SSlsysm}, with a right hand side, $F$, such that
\begin{equation}\label{FZest}
\big(F,Z[\underline U] U\big) \ \leq \ C_F \ \max(\epsilon,\beta) \big\vert U \big\vert_{X^0}^2\ +\ f(t)\ \big\vert U \big\vert_{X^0}\nn,\end{equation}
with $C_F$ a constant and $f$ a positive integrable function on $[0,T/\max(\epsilon,\beta)]$.\\
Then there exists $\lambda,C_1\equiv C(\big\Vert \partial_t \underline{U} \big\Vert_{L^\infty},\big\Vert \underline U \big\Vert_{L^\infty},\big\Vert \partial_x\underline U \big\Vert_{L^\infty},\Vert b\Vert_{H^{s_0+3}},C_F)$ such that
\begin{equation}\label{energyestimateL2Z}
	\forall t\in [0,\frac{T}{\max(\epsilon,\beta)}],\qquad
	E^0(U)(t)\leq e^{\max(\epsilon,\beta)\lambda t}E^0(U_0)+ \int^{t}_{0} e^{\max(\epsilon,\beta)\lambda( t-t')}\Big(f(t')+\max(\epsilon,\beta) C_1 \Big) dt',
\end{equation}
The constants $\lambda$ and $C_1$ are independent of $\p=(\mu, \epsilon, \delta, \gamma, \beta, \bo)\in \P_{\rm CH}$, but depend on $M_{\rm CH},h_{01}^{-1},h_{02}^{-1}$, and $h_{03}^{-1}$.
\end{lemma}
\emph{Proof}.\\
Let us take the inner product of~\eqref{SSlsysm} by $ Z[\underline U] U$:
\[ \big(\partial_t U,Z[\underline U] U\big) \ + \ \big(A[\underline U]\partial_x U,Z[\underline U]  U\big)\ + \ \big(B[\underline U],Z[\underline U]  U\big)
\ =  \ \big( F ,Z[\underline U] U\big) \ ,
\]
From the symmetry property of $Z[\underline U]$, and using the  definition of $E^s(U)$, one deduces
\begin{multline}\label{qtsctrL2Z}
\frac12 \frac{d}{dt}E^0(U) ^2  =  \frac12\big( U,\big[\partial_t,  Z[\underline U]\big] U\big)-\big(Z[\underline U]A[\underline U]\partial_xU,  U\big)
-\big(B[\underline U],Z[\underline U]  U\big)+\big(F,Z[\underline U] U\big).
\end{multline}

Let us first estimate $\big(B[\underline U],Z[\underline U]  U\big)$. One has
\begin{eqnarray*}
\big(B[\underline U],Z[\underline U]  U\big) &  = &  \Big( - g(\epsilon\underline{\zeta})\underline{v}\beta\partial_x b,\frac{Q_0(\epsilon\underline{\zeta},\beta b)+\epsilon^2Q_1(\epsilon\underline{\zeta},\beta b,\underline{v})}{f(\epsilon\underline{\zeta})}\zeta\Big)\\ &\qquad+&\Big(\epsilon\mfT^{-1}\big(\dfrac{\gamma\beta q_1(\epsilon\zetabar,\beta b)h_1(h_1+h_2)\underline{v}^2\partial_x b}{(h_1+\gamma h_2)^3}\big),\mfT[\epsilon\zetabar,\beta b]v\Big).
 \end{eqnarray*}

  \[
 \Big( - g(\epsilon\underline{\zeta})\underline{v}\beta\partial_x b,\frac{Q_0(\epsilon\underline{\zeta},\beta b)+\epsilon^2Q_1(\epsilon\underline{\zeta},\beta b,\underline{v})}{f(\epsilon\underline{\zeta})}\zeta\Big) \leq \beta C\big(\big\Vert \underline{U} \big\Vert_{L^\infty},\Vert b\Vert_{W^{2,\infty}},\Vert \partial_x b\Vert_{L^2}\big)\big\vert U\big\vert_{X^0}.
\]
Using the symmetry property of $\mfT[\epsilon\zetabar,\beta b]$, we write
$$\Big(\epsilon\mfT^{-1}\big(\dfrac{\gamma\beta q_1(\epsilon\zetabar,\beta b)h_1(h_1+h_2)\underline{v}^2\partial_x b}{(h_1+\gamma h_2)^3}\big),\mfT[\epsilon\zetabar,\beta b]v\Big)=
\epsilon\Big(\dfrac{\gamma\beta q_1(\epsilon\zetabar,\beta b)h_1(h_1+h_2)\underline{v}^2\partial_x b}{(h_1+\gamma h_2)^3},v\Big).$$
From Cauchy-Schwarz inequality, one deduces
$$\epsilon\Big(\dfrac{\gamma\beta q_1(\epsilon\zetabar,\beta b)h_1(h_1+h_2)\underline{v}^2\partial_x b}{(h_1+\gamma h_2)^3},v\Big)\leq  \beta C\big(\big\Vert \underline{U} \big\Vert_{L^\infty},\Vert \partial_x b\Vert_{L^2}\big)\big\vert U\big\vert_{X^0}.$$

Altogether, one has
\begin{equation}\label{BZest}
\big(B[\underline U],Z[\underline U]  U\big) \leq\beta C_1 \big\vert U\big\vert_{X^0}\leq \max(\epsilon,\beta) C_1 \big\vert U\big\vert_{X^0}.
\end{equation}
Now we have,
 \begin{equation*}
Z[\underline{U}]A[\underline{U}]=\begin{pmatrix}
 \epsilon\frac{Q_0(\epsilon\underline{\zeta},\beta b)+\epsilon^2Q_1(\epsilon\underline{\zeta},\beta b,\underline{v}) }{f(\epsilon\underline{\zeta})}f'(\epsilon\underline{\zeta}) \underline{v}& \quad Q_0(\epsilon\underline{\zeta},\beta b)+\epsilon^2Q_1(\epsilon\underline{\zeta},\beta b,\underline{v}) \\ \\
Q_0(\epsilon\underline{\zeta},\beta b)+\epsilon^2 Q_1(\epsilon\underline{\zeta},\beta b,\underline{v})  &\quad \epsilon \mfQ[\epsilon\underline{\zeta},\beta b,\underline{v}]+\epsilon \varsigma \mfT[\epsilon\underline{\zeta},\beta b](\underline{v} .)
\end{pmatrix}
\end{equation*}
One has,
\begin{eqnarray*}
\big(Z[\underline{U}]A[\underline{U}]\partial_x U, U\big)&=&\Big(\epsilon\dfrac{Q_0(\epsilon\zetabar,\beta b)+\epsilon^2Q_1(\epsilon\zetabar,\beta b,\vbar)}{f(\epsilon\zetabar)}f'(\epsilon\zetabar)\vbar\partial_x \zeta, \zeta \Big)\\&+& \Big(Q_0(\epsilon\zetabar,\beta b)\partial_x v,\zeta\Big)+\Big(\epsilon^2Q_1(\epsilon\zetabar,\beta b,\vbar)\partial_x v,\zeta \Big)\\&+& \Big(Q_0(\epsilon\zetabar,\beta b)\partial_x \zeta, v \Big) +\Big(\epsilon^2Q_1(\epsilon\zetabar,\beta b,\vbar)\partial_x \zeta, v \Big)\\&+& \Big(\epsilon \mfQ[\epsilon\underline{\zeta},\beta b,\underline{v}] \partial_x v, v \Big)+ \epsilon \varsigma\Big( \mfT[\epsilon\underline{\zeta},\beta b](\underline{v}\partial_x v), v \Big).
\end{eqnarray*}
One deduces that,
\begin{eqnarray*}
\big(Z[\underline{U}]A[\underline{U}]\partial_x U, U\big)
&=&-\dfrac{1}{2}\Big(\epsilon\partial_x\big(\dfrac{Q_0(\epsilon\zetabar,\beta b)+\epsilon^2Q_1(\epsilon\zetabar,\beta b,\vbar)}{f(\epsilon\zetabar)}f'(\epsilon\zetabar)\vbar\big)\zeta, \zeta \Big)\\&-&\Big(\partial_x\big(Q_0(\epsilon\zetabar,\beta b)\big) \zeta ,v \Big)-\epsilon^2\Big(\partial_x\big(Q_1(\epsilon\zetabar,\beta b,\vbar)\big)\zeta, v\Big)\\&+&\Big(\epsilon \mfQ[\epsilon\underline{\zeta},\beta b,\underline{v}] \partial_x v, v \Big)+ \epsilon \varsigma\Big( \mfT[\epsilon\underline{\zeta},\beta b](\underline{v}\partial_x v), v \Big).
\end{eqnarray*}
One can easily remark that we didn't use the smallness assumption of the Camassa-Holm regime $\epsilon=\OO(\sqrt\mu)$ since we do not have anymore $\partial_xv$ in the third term of the above identity.\\
\\
One make use of the identity below,

 \begin{equation*}\label{eqn:symmetry-like} \begin{array}{rcl}\big(\mfT[\epsilon\underline{\zeta},\beta b] (\underline{v} \partial_x V) ,V\big)&=& \big(\ q_1(\epsilon\underline{\zeta},\beta b)\underline{v} \partial_x V-\mu\nu\partial_x(q_2(\epsilon\underline{\zeta},\beta b)\partial_x(\underline{v} \partial_x V ))\ ,\ V\ \big)  \\
&=& -\frac12\big(\ \partial_x(q_1(\epsilon\underline{\zeta},\beta b)\underline{v}) V\ ,\ V\ \big)+\mu\nu\big(\ q_2(\epsilon\underline{\zeta},\beta b)\partial_x(\underline{v} \partial_x V )\ ,\ \partial_x V\ \big)\\
&=& -\frac12\big( \partial_x(q_1(\epsilon\underline{\zeta},\beta b)\underline{v}) V , V \big)+\mu\nu\big( q_2(\epsilon\underline{\zeta},\beta b)(\partial_x\underline{v}) \partial_x V  , \partial_x V \big)\\&&\qquad -\mu\nu\frac12\big( \partial_x( q_2(\epsilon\underline{\zeta},\beta b)\underline{v} )\partial_x V  , \partial_x V \big).
\end{array}\end{equation*}

From Cauchy-Schwarz inequality, one deduces
\begin{equation}\label{estZA0A1}
\Big|\big(Z[\underline{U}]A[\underline{U}]\partial_x U, U\big)|\leq\max(\epsilon,\beta) C\Big(\big\Vert \underline U \big\Vert_{L^\infty}+\big\Vert \partial_x \underline U \big\Vert_{L^\infty}+\big\Vert b \big\Vert_{W^{3,\infty}}\Big) \big\vert U\big\vert_{X^0}^2.
\end{equation}
The last term to estimate is $\big(  U,\big[\partial_t,  Z[\underline U]\big]  U\big)$.\\
\\
One has
\begin{eqnarray*}
\big(  U,\big[\partial_t,  Z[\underline U]\big]  U\big)&\equiv&\dsp( v,\Big[\partial_t,\underline{\mfT}\Big] v)
+( \zeta,\Big[\partial_t,\frac{Q_0(\epsilon\underline{\zeta},\beta b)+\epsilon^2Q_1(\epsilon\zetabar,\beta b,\vbar)}{f(\epsilon\underline{\zeta})}\Big] \zeta)\\
&=&\Big(  v,\big(\partial_t q_1(\epsilon\underline{\zeta},\beta b) \big)  v\Big)-\mu\nu\Big(  v,\partial_x\big((\partial_tq_2(\epsilon\underline{\zeta},\beta b))( \partial_x  v)\big) \Big)\\
&&\qquad \qquad+\Big( \zeta,\partial_t \Big(\frac{Q_0(\epsilon\underline{\zeta},\beta b)+\epsilon^2 Q_1(\epsilon\zetabar,\beta b,\vbar)}{f(\epsilon\underline{\zeta})}\Big) \zeta\Big).\\
\end{eqnarray*}

From Cauchy-Schwarz inequality and since $\underline{\zeta}$ and $b$ satisfies~\eqref{CondDepth}, one deduces
\begin{eqnarray}\label{partialZest}
\big\vert \ \frac12\big(  U,\big[\partial_t, Z[\underline U]\big]  U\big)\ \big\vert\ &\leq& \  \epsilon C(\big\Vert \partial_t \underline U \big\Vert_{L^\infty},\big\Vert \underline U \big\Vert_{L^\infty})\big\vert U\big\vert_{X^0}^2 \nn \\  &\leq& \  \max(\epsilon,\beta) C(\big\Vert \partial_t \underline U \big\Vert_{L^\infty},\big\Vert \underline U \big\Vert_{L^\infty})\big\vert U\big\vert_{X^0}^2. \nn \\
\end{eqnarray}
One can now conclude with the proof of the $X^0$ energy estimate. Plugging \eqref{BZest}, \eqref{estZA0A1} and \eqref{partialZest} into \eqref{qtsctrL2Z}, and making use of the assumption of the Lemma on F. This yields
\[ \frac12 \frac{d}{dt}E^0(U)^2 \ \leq \ \max(\epsilon,\beta) \ C_1 E^0(U)^2+ \Big(f(t)+\max(\epsilon,\beta) C_1 \Big) E^0(U),\]
where $C_1 \equiv C(\big\Vert \partial_t \underline{U} \big\Vert_{L^\infty},\big\Vert \underline U \big\Vert_{L^\infty},\big\Vert \partial_x\underline U \big\Vert_{L^\infty},\Vert b\Vert_{H^{s_0+3}},C_F)$ . Consequently
\[ \frac{d}{dt}E^0(U) \ \leq \ \max(\epsilon,\beta)C_1E^0(U)+ \Big(f(t)+\max(\epsilon,\beta) C_1 \Big).\]
Making use of the usual trick, we compute for any $\lambda \in \RR$,
\[
e^{\max(\epsilon,\beta)\lambda t}\partial_t(e^{-\max(\epsilon,\beta)\lambda t}E^0(U))=-\max(\epsilon,\beta)\lambda E^0(U) +\frac{d}{dt} E^0(U).
\]
Thanks to the above inequality, one can choose  $\lambda= C_1$, so that for all $t\in [0,\frac{T}{\max(\epsilon,\beta)}]$, one deduces
\[
\frac{d}{dt} (e^{-\max(\epsilon,\beta)\lambda t}E^0(U)) \ \leq \  \Big(f(t)+\max(\epsilon,\beta) C_1 \Big)e^{-\max(\epsilon,\beta)\lambda t} .\]
Integrating this differential inequality yields
\begin{equation}\label{energyestimateL2inproof}
	\forall t\in [0,\frac{T}{\max(\epsilon,\beta)}],\qquad
	E^0(U)(t)\leq e^{\max(\epsilon,\beta)\lambda t}E^0(U_0)+ \int^{t}_{0} e^{\max(\epsilon,\beta)\lambda( t-t')}\Big(f(t')+\max(\epsilon,\beta) C_1 \Big)dt'.
\end{equation}
This proves the energy estimate~\eqref{energyestimateL2Z}.\demo\\
\\
Let us now turn to the {\em a priori} energy estimate in ``large'' $X^s$ norm.
\begin{lemma}[ $X^s$ energy estimate]\label{Lem:HsZ}
Set $(\mu, \epsilon, \delta, \gamma, \beta, \bo)\in \P_{\rm CH}$, and  $s\geq s_0+1$, $s_0>1/2$. Let $U=(\zeta,v)^\top$ and $\underline{U}=(\underline{\zeta},\underline{v})^\top$ be such that $U,\underline{U}\in L^\infty([0,T/\max(\epsilon,\beta)];X^{s}) $, $\partial_t\underline{U}\in L^\infty([0,T/\max(\epsilon,\beta)]\times\RR)$, $b \in H^{s+2}$ and $\underline{U}$ satisfies~\eqref{CondDepth},\eqref{CondEllipticity}, and \eqref{H3} uniformly on $[0,T/\max(\epsilon,\beta)]$, and such that system~\eqref{SSlsysm} holds with a right hand side, $F$, with
\[
\big( \Lambda^s F,Z[\underline U]\Lambda^s U\big)\ \leq  \ C_F\ \max(\epsilon,\beta) \big\vert  U\big\vert_{X^s}^2\ +\ f(t)\ \big\vert U\big\vert_{X^s} ,
\]
where $C_F$ is a constant and $f$ is an integrable function on $[0,T/\max(\epsilon,\beta)]$.\\
\\
Then there exists $\lambda,C_2=C(\big\Vert \underline U\big\Vert_{X^s_T},\big\Vert b \big\Vert_{H^{s+2}},C_F)$ such that the following energy estimate holds:
\begin{equation}\label{energyestimateZs}
E^s(U)(t)\leq e^{\max(\epsilon,\beta)\lambda t}E^s(U_0)+
\int^{t}_{0} e^{\max(\epsilon,\beta)\lambda( t-t')}\big(f(t')+\max(\epsilon,\beta)C_2\big) dt',\end{equation}

The constants $\lambda$ and $C_2$ are independent of $\p=(\mu, \epsilon, \delta, \gamma, \beta, \bo)\in \P_{\rm CH}$, but depend on $M_{\rm CH},h_{01}^{-1},h_{02}^{-1}$, and $h_{03}^{-1}$.
\end{lemma}
\begin{remark}
In this Lemma, and in the proof below, the norm $\big\Vert \underline U\big\Vert_{X^s_T}$ is to be understood as essential sup:
\[
 \Vert U\Vert_{X^s_T}\equiv \esssup_{t\in [0,T/\max(\epsilon,\beta)]}\vert U(t,\cdot)\vert_{X^s}+\esssup_{t\in [0,T/\max(\epsilon,\beta)],x\in\RR}\vert \partial_t U(t,x)\vert.
\]
\end{remark}
\emph{Proof}.\\
Let us multiply the system~\eqref{SSlsysm} on the right by $\Lambda^s Z[\underline U]\Lambda^s U$, and integrate by parts. One obtains
\begin{multline} \big(\Lambda^s\partial_t U,Z[\underline U]\Lambda^s U\big) \ + \ \big(\Lambda^sA[\underline U]\partial_x U,Z[\underline U]\Lambda^s U\big)\ + \  \big(\Lambda^s B[\underline U],Z[\underline U]\Lambda^s U\big) 
\ =  \ \big(\Lambda^s F,Z[\underline U]\Lambda^s U\big)\ ,
\end{multline}
from which we deduce, using the symmetry property of $Z[\underline U]$, as well as the definition of $E^s(U)$:
\begin{multline}\label{qtsctrsZ}
\frac12 \frac{d}{dt}E^s(U)^2   \ =  \frac12\big(\Lambda^s U,\big[\partial_t,  Z[\underline U]\big]\Lambda^s U\big)-\big(Z[\underline U]A[\underline U]\partial_x \Lambda^s U, \Lambda^s U\big) 
-\big(\big[\Lambda^s ,A[\underline U]\big]\partial_x U,Z[\underline U]\Lambda^s U\big)\\
-\big(\Lambda^s B[\underline U],Z[\underline U]\Lambda^s U\big)+\big( \Lambda^s F,Z[\underline U]\Lambda^s U\big).
\end{multline}
We now estimate each of the different components of the r.h.s of the above identity.

\bigskip

\noindent $\bullet$ {\em Estimate of  $\big(\Lambda^s B[\underline U],Z[\underline U] \Lambda^s  U\big)$,}
\begin{eqnarray*}
\big(\Lambda^s B[\underline U],Z[\underline U] \Lambda^s  U\big)= \Big( \Lambda^s\big(- g(\epsilon\underline{\zeta})\underline{v}\beta\partial_x b\big),\frac{Q_0(\epsilon\underline{\zeta},\beta b)+\epsilon^2 Q_1(\epsilon\underline{\zeta},\beta b,\underline{v}) }{f(\epsilon\underline{\zeta})}\Lambda^s\zeta\Big) \\+\epsilon\Big(\Lambda^s\underline{\mfT}^{-1}\big(\dfrac{\gamma\beta q_1(\epsilon\zetabar,\beta b)h_1(h_1+h_2)\underline{v}^2}{(h_1+\gamma h_2)^3}\partial_x b\big),\mfT[\epsilon\zetabar,\beta b]\Lambda^s v\Big).
\end{eqnarray*}
Using Cauchy-Schwarz inequality, Lemma~4.2 and Lemma~4.7 one has,
\begin{eqnarray}\label{estBZs}|\big(\Lambda^s B[\underline U],Z[\underline U] \Lambda^s  U\big)|\leq \beta C\big(\big\Vert \underline{U} \big\Vert_{X^s_T},\Vert b\Vert_{H^{s+1}}\big)|U| _{X^s}\leq\max(\epsilon,\beta) C_2|U| _{X^s}\nn.\\\end{eqnarray}

\noindent $\bullet$ {\em Estimate of $\big(Z[\underline U]A[\underline U]\partial_x \Lambda^s U, \Lambda^s U\big) $.}\\
\\
Thanks to Sobolev embedding, one has for $s>s_0+1, s_0>1/2$
\[C(\Vert\underline{U}\Vert_{L^{\infty}}+\Vert\partial_x\underline{U}\Vert_{L^{\infty}})\leq C(\big\Vert \underline U \big\Vert_{X^s_T}),\]
One can use the $L^2$ estimate derived in~\eqref{estZA0A1}, applied to $\Lambda^s U$. One deduces
\begin{equation}\label{estZA0A1s}
\Big|\big(Z[\underline{U}]A[\underline{U}]\partial_x\Lambda^s U, \Lambda^sU\big)\Big|\leq\max(\epsilon,\beta) C\Big(\big\Vert \underline U \big\Vert_{X^s_T}+\big\Vert b \big\Vert_{W^{3,\infty}}\Big) \big\vert U\big\vert_{X^s}^2.
\end{equation}

\noindent $\bullet$ {\em Estimate of $ \big(\big[\Lambda^s,A[\underline{U}]\big]\partial_x U,Z[\underline{U}]\Lambda^s U\big)$}. Using the definition of $A[\cdot]$ and $Z[\cdot]$ in~\eqref{defA0A1} and \eqref{defZ}, one has
 \begin{eqnarray*}
\big(\big[\Lambda^s,A[\underline{U}]\big]\partial_x U,Z[\underline{U}]\Lambda^s U\big)=\Big([\Lambda^s,\epsilon f'(\epsilon\underline{\zeta}) \underline{v}]\partial_x\zeta +[\Lambda^s, f(\epsilon\underline{\zeta})]\partial_x v\ ,\ \frac{Q(\epsilon\underline{\zeta},\beta b,\underline{v})}{f(\epsilon\underline{\zeta})}\Lambda^s\zeta\Big) \\
+\Big([\Lambda^s, \underline{\mfT}^{-1}\big(Q(\epsilon\underline{\zeta},\beta b,\underline{v})\big)]\partial_x \zeta\ ,\ \underline{\mfT}\Lambda^s v\Big)\  +\ \epsilon\Big([\Lambda^s, \underline{\mfT}^{-1}\mfQ[\epsilon\underline{\zeta},\beta b,\underline{v}] +\varsigma \underline v]\partial_x v, \underline{\mfT}\Lambda^s v\Big).
\end{eqnarray*}
Here and in the following, we denote  $\underline{\mfT}\equiv \mfT[\epsilon\underline{\zeta},\beta b]$ and $Q(\epsilon\underline{\zeta},\beta b,\underline{v})=Q_0(\epsilon\underline{\zeta},\beta b)+\epsilon^2 Q_1(\epsilon\underline{\zeta},\beta b,\underline{v}).$\\
\\
Using the same techniques as in~\cite[Lemma 6.6 ]{DucheneIsrawiTalhouk13} and since $\underline{\zeta}$ and $b$ satisfies \eqref{CondDepth}, we proved
\begin{eqnarray}\label{estAZs}
\Big|\big(\big[\Lambda^s,A[\underline{U}]\big]\partial_x U,Z[\underline{U}]\Lambda^s U\big)\Big|&\leq& \max(\epsilon, \beta) C\Big(\big\Vert \underline U \big\Vert_{X^s_T}+\big\Vert b \big\Vert_{H^{s+2}} \Big)\big\vert U\big\vert_{X^s}^2.\nn \\
\end{eqnarray}

\noindent $\bullet$ {\em Estimate of $\frac12\big(\Lambda^s U,\big[\partial_t,  Z[\underline U]\big]\Lambda^s U\big)$.}

One has
\begin{eqnarray*}
\big(\Lambda^s U,\big[\partial_t,  Z[\underline U]\big]\Lambda^s U\big)&\equiv&\dsp(\Lambda^sv,\Big[\partial_t,\underline{\mfT}\Big]\Lambda^sv)
+(\Lambda^s\zeta,\Big[\partial_t,\frac{Q_0(\epsilon\underline{\zeta},\beta b)+\epsilon^2 Q_1(\epsilon\underline{\zeta},\beta b,\underline{v})}{f(\epsilon\underline{\zeta})}\Big]\Lambda^s\zeta)\\
&=&\Big(\Lambda^s v,\big(\partial_t q_1(\epsilon\underline{\zeta},\beta b) \big)\Lambda^s v\Big)-\mu\nu\Big(\Lambda^s v,\partial_x\big((\partial_tq_2(\epsilon\underline{\zeta},\beta b))( \partial_x\Lambda^s v)\big) \Big)\\
&&\qquad \qquad
+\Big(\Lambda^s\zeta,\partial_t \Big(\frac{Q_0(\epsilon\underline{\zeta},\beta b)+\epsilon^2 Q_1(\epsilon\underline{\zeta},\beta b,\underline{v})}{f(\epsilon\underline{\zeta})}\Big)\Lambda^s\zeta\Big)\\
&=&\epsilon\kappa_1\Big(\Lambda^s v,(\partial_t \underline{\zeta}) \Lambda^s v\Big)+\mu\nu\epsilon\kappa_2\Big(\Lambda^s \partial_x v,(\partial_t\underline{\zeta})\Lambda^s \partial_x v \Big)\\
&&\qquad +\Big(\Lambda^s\zeta,\partial_t \Big(\frac{Q_0(\epsilon\underline{\zeta},\beta b)+\epsilon^2 Q_1(\epsilon\underline{\zeta},\beta b,\underline{v})}{f(\epsilon\underline{\zeta})}\Big)\Lambda^s\zeta\Big).\end{eqnarray*}

From Cauchy-Schwarz inequality and since $\underline{\zeta}$ and $b$ satisfies~\eqref{CondDepth}, one deduces
\begin{eqnarray}
\Big|\ \frac12\big( \Lambda^s U,\big[\partial_t, Z[\underline U]\big] \Lambda^s U\big)\ \Big| &\leq& \  \epsilon C(\big\Vert \partial_t \underline U \big\Vert_{L^\infty},\big\Vert \underline U \big\Vert_{L^\infty})\big\vert U\big\vert_{X^s}^2 \nn \\  &\leq& \  \max(\epsilon,\beta) C(\big\Vert \partial_t \underline U \big\Vert_{L^\infty},\big\Vert \underline U \big\Vert_{L^\infty})\big\vert U\big\vert_{X^s}^2. \nn
\end{eqnarray}
and continuous Sobolev embedding yields,
\begin{eqnarray}\label{partialZests}
\Big|\ \frac12\big( \Lambda^s U,\big[\partial_t, Z[\underline U]\big] \Lambda^s U\big)\ \Big| \leq \  \max(\epsilon,\beta) C\big(\big\Vert \underline U\big\Vert_{X^s_T}\big)\big\vert U\big\vert_{X^s}^2.
\end{eqnarray}
One can now conclude the proof of the $X^s$ energy estimate. Plugging \eqref{estBZs}, \eqref{estZA0A1s}, \eqref{estAZs} and \eqref{partialZests} into \eqref{qtsctrsZ}, and making use of the assumption of the Lemma on F.
\[ \frac12 \frac{d}{dt}E^s(U)^2 \ \leq \  \max(\epsilon,\beta)C_2  E^s(U)^2+ E^s(U)\big(f(t)+\max(\epsilon,\beta)  C_2 \big),\]
with $C_2=C(\big\Vert \underline U\big\Vert_{X^s_T},\big\Vert b \big\Vert_{H^{s+2}},C_F)$, and consequently
\[ \frac{d}{dt}E^s(U) \ \leq \ \max(\epsilon,\beta) C_2 E^s(U)+ \big(f(t)+\max(\epsilon,\beta)  C_2 \big) .\]
Making use of the usual trick, we compute for any $\lambda \in \RR$,
\[
e^{\max(\epsilon,\beta) \lambda t}\partial_t(e^{-\max(\epsilon,\beta) \lambda t}E^s(U))=-\max(\epsilon,\beta) \lambda E^s(U) +\frac{d}{dt} E^s(U).
\]
Thus with $\lambda=C_2$, one has for all $t\in [0,\frac{T}{\max(\epsilon,\beta) }]$,
\[
\frac{d}{dt} (e^{-\max(\epsilon,\beta) \lambda t}E^s(U)) \ \leq \  \big(f(t)+\max(\epsilon,\beta) C_2 \big )e^{-\max(\epsilon,\beta) \lambda t} .\]
Integrating this differential inequality yields,
\[E^s(U)(t)\leq e^{\max(\epsilon,\beta) \lambda t}E^s(U_0)+
\int^{t}_{0} e^{\max(\epsilon,\beta) \lambda( t-t')}\big(f(t')+\max(\epsilon,\beta)  C_2 \big)dt'.\]\demo

\subsection{Well-posedness of the linearized system}

\begin{prop}\label{ESpropZ}
Let $\p=(\mu, \epsilon, \delta, \gamma, \beta, \bo)\in \P_{\rm CH}$ and $s\geq s_0+1$ with $s_0>1/2$, and let $\underline{U}=(\underline{\zeta}, \underline{v})^\top\in X^{s}_T$ (see Definition~\ref{defispace}), $b \in H^{s+2}$ be such that~\eqref{CondDepth},\eqref{CondEllipticity}, and \eqref{H3} are satisfied for $t\in [0,T/\max(\epsilon,\beta)]$, uniformly with respect to $\p \in \P_{\rm CH}$. For any $U_0\in X^{s}$, there exists
a unique solution to~\eqref{SSlsys}, $U^{\p} \in C^0([0,T/\max(\epsilon,\beta)]; X^{s})\cap C^1([0,T/\max(\epsilon,\beta)]; X^{s-1})\subset X^s_T$,\\  with
$\lambda_T,C_0= C(\big\Vert \underline U\big\Vert_{ X^{s}_T},T,M_{\rm CH},h_{01}^{-1},h_{02}^{-1},h_{03}^{-1},\Vert b \Vert_{H^{s+2}})$, independent of $\p\in \P_{\rm CH}$, such that
 one has the energy estimates
\[
\forall\  0\leq t\leq\frac{T}{\max(\epsilon,\beta)}, \quad E^s(U^{\p})(t)\leq e^{\max(\epsilon,\beta)\lambda_{T} t}E^s(U_0)+\max(\epsilon,\beta)C_0 \int_{0}^{t}e^{\max(\epsilon,\beta)\lambda_T(t-t')}dt'
\]
\[\mbox{and} \ \ \ E^{s-1}( \partial_t U^{\p}) \leq C_0 e^{\max(\epsilon,\beta)\lambda_T t}E^s(U_0)+\max(\epsilon,\beta)C_0^2 \int_{0}^{t} e^{\max(\epsilon,\beta)\lambda_T (t-t')}dt'+\max(\epsilon,\beta)  C_0. \]
\end{prop}
\emph{Proof}.\\
Existence and uniqueness of a solution to the initial value problem~\eqref{SSlsys} follows, by standard techniques, from the estimate~\eqref{energyestimateZs} in Lemma~\ref{Lem:HsZ}:
\begin{equation}\label{energyestimateinproof}
E^s(U)(t)\leq e^{\max(\epsilon,\beta)\lambda_T t}E^s(U_0)+\max(\epsilon,\beta) C_0 \int_{0}^{t} e^{\max(\epsilon,\beta)\lambda_T (t-t')}dt',\end{equation}
(since $F\equiv0$, and omitting the index $\p$ for the sake of simplicity).\\
First, let us notice that using the system of equation~\eqref{SSlsys}, one can deduce an energy estimate on the time-derivative of the solution. Indeed, one has
\begin{eqnarray}
\big\vert\partial_t U\big\vert_{X^{s-1}}&=& \big\vert -A[\underline{U}]\partial_x U-B[\underline{U}]\big\vert_{X^{s-1}}\nn \\
&\leq &\big\vert \epsilon f'(\epsilon\underline{\zeta})\underline{v}\partial_x \zeta +f(\epsilon\underline{\zeta})\partial_x v+ \beta \partial_x b g(\epsilon\zetabar) \underline{v}\big\vert_{H^{s-1}}\nn \\
&&\quad + \big\vert  \mathfrak{T}[\epsilon\underline{\zeta},\beta b]^{-1}\Big(Q_0(\epsilon\underline{\zeta},\beta b)\partial_x\zeta+\epsilon \mathfrak{Q}[\epsilon\underline{\zeta},\beta b,\underline{v}]\partial_x v+\epsilon^2Q_1(\epsilon\underline{\zeta},\beta b,\underline{v})\partial_x\zeta \nn \\ &&\quad +\epsilon \dfrac{\gamma\beta q_1(\epsilon\zetabar,\beta b)h_1(h_1+h_2)\underline{v}^2\partial_x b}{(h_1+\gamma h_2)^3} \Big)+\epsilon\varsigma\underline{v}\partial_x v\big\vert_{H^s_\mu} \nn \\
&\leq& C(| \underline{U}|_{X^s}, |b|_{H^{s+1}})|U|_{X^s}+\beta C_0\nn\\
&\leq& C_0 E^s(U)(t)+\beta C_0\nn\\
&\leq& C_0 e^{\max(\epsilon,\beta)\lambda_{T} t}E^s(U_0)+\max(\epsilon,\beta)C_0^2 \int_{0}^{t} e^{\max(\epsilon,\beta)\lambda_T (t-t')}dt'+\max(\epsilon,\beta)C_0.\nn\\
\label{energyestimatederivativeinproof}\end{eqnarray}
The completion of the proof is as follows. In order to construct a solution to~\eqref{SSlsys}, we use a sequence of Friedrichs mollifiers, defined by $J_\nu\equiv (1-\nu\partial_x^2)^{-1/2}$ ($\nu>0$), in order to reduce our system to ordinary differential equation systems on $X^s$, which are solved uniquely by Cauchy-Lipschitz theorem. Estimates~\eqref{energyestimateinproof},\eqref{energyestimatederivativeinproof} hold for each $U_\nu\in C^0([0,T/\max(\epsilon,\beta)];X^{s})$, uniformly in $\nu>0$. One deduces that a subsequence converges towards $U\in L^2([0,T/\max(\epsilon,\beta)];X^s)$, a (weak) solution of the Cauchy problem~\eqref{SSlsys}. By regularizing the initial data as well, one can show that the system induces a smoothing effect in time, and that the solution $U\in C^0([0,T/\max(\epsilon,\beta)];X^{s})\cap C^1([0,T/\max(\epsilon,\beta)];X^{s-1})$ is actually a strong solution. The uniqueness is a straightforward consequence of~\eqref{energyestimateinproof} (with $U_0\equiv 0$) applied to the difference of two solutions.\demo

\subsection{A priori estimate}
In this section, we control the difference of two solutions of the nonlinear system, with different initial data and right-hand sides.
\begin{prop} \label{prop:stabilityZ}
 Let $(\mu, \epsilon, \delta, \gamma, \beta, \bo)\in \P_{\rm CH}$ and $s\geq s_0+1$, $s_0>1/2$, and assume that there exists $U_i$ for $i\in \{1,2\}$, such that
 $U_{i}=(\zeta_{i},v_{i})^\top \in  X^{s}_{T}$, $U_2\in L^\infty([0,T/\max(\epsilon,\beta) ];X^{s+1})$, $b\in H^{s+2}$, $U_1$ satisfy~\eqref{CondDepth},\eqref{CondEllipticity} and \eqref{H3} on $[0,T/\max(\epsilon,\beta) ]$, with $h_{01},h_{02},h_{03}>0$,
 and $U_i$ satisfy
\begin{eqnarray*}
		\partial_t U_1\ +\ A[U_1]\partial_x U_1\ +\ B[U_1] &=&\ F_1 \ ,  \\
		\partial_t U_2\ +\ A[U_2]\partial_x U_2\ +\ B[U_2]  &=&\ F_2 \ ,
\end{eqnarray*}		
with $F_i\in L^1([0,T/\max(\epsilon,\beta) ];X^{s})$.\\
\\
Then there exists constants $C_0=C(M_{\rm CH},h_{01}^{-1},h_{02}^{-1},h_{03}^{-1},\max(\epsilon,\beta) \big\vert U_1\big\vert_{X^s},\max(\epsilon,\beta) \big\vert U_2\big\vert_{X^s},|b|_{H^{s+2}})$ \\ and $\lambda_T=\big(C_0\times C( \vert U_2\vert_{L^\infty([0,T/\max(\epsilon,\beta) ];X^{s+1})})+C_0\big)$  such that for all $t \in [0,\frac{T}{\max(\epsilon,\beta) }]$,
\begin{eqnarray*}
E^s(U_1-U_2)(t)\leq e^{\max(\epsilon,\beta) \lambda_{T} t}E^s(U_{1}\id{t=0}-U_{2}\id{t=0})+C_0 \int^{t}_{0} e^{\max(\epsilon,\beta) \lambda_T( t-t')} E^s(F_1-F_2)(t')dt'.
\end{eqnarray*}
\end{prop}
\emph{Proof}.\\
When multiplying the equations satisfied by $U_i$ on the left by $Z[U_i]$, one obtains
	\begin{eqnarray*}
	Z[U_1]\dsp\partial_t U_1+\Sigma[U_1]\partial_x U_1+Z[U_1]B[U_1]&=& Z[U_1]F_1\\
	Z[U_2]\dsp\partial_t U_2+\Sigma[U_2]\partial_x U_2
+Z[U_2]B[U_2]&=& Z[U_2]F_2;
    \end{eqnarray*}
with $\Sigma[U]=Z[U]A[U]$.
Subtracting the two equations above, and defining  $V=U_1-U_2\equiv (\zeta,v)^\top$ one obtains
\begin{multline*}
	Z[U_1]\dsp\partial_t V+\Sigma[U_1]\partial_x V+(Z[U_1]B[U_1]-Z[U_2]B[U_2])=
	 Z[U_1](F_1-F_2)-(\Sigma[U_1]-\Sigma[U_2])\partial_xU_2\\-(Z[U_1]-Z[U_2])(\partial_tU_2-F_2).
\end{multline*}
We then apply $Z^{-1}[U_1]$ and deduce the following system satisfied by $V$:
\begin{equation}\label{systemstabilityZ}
	\left\lbrace
	\begin{array}{l}
	\dsp\partial_t V+A[U_1]\partial_x V+Z^{-1}[U_1]\big(Z[U_1]B[U_1]-Z[U_2]B[U_2]\big)=F \\
	\dsp V(0)=(U_{1}-U_{2})\id{t=0},
	\end{array}\right.
\end{equation}
\begin{equation}\label{rhsstabilityZ}
\mbox{where,}\quad F \equiv F_1-F_2 - Z^{-1}[U_1]\big(\Sigma[U_1]-\Sigma[U_2]\big)
\partial_xU_2- Z^{-1}[U_1]\big(Z[U_1]-Z[U_2]\big)(\partial_tU_2-F_2).
\end{equation}
We wish to use the energy estimate of Lemma~\ref{Lem:HsZ} to the linear system~\eqref{systemstabilityZ}.\\
\\
The additional term now is $Z^{-1}[U_1]\big(Z[U_1]B[U_1]-Z[U_2]B[U_2]\big)$.\\
\\
So we have to control, $$\big(\Lambda^s Z^{-1}[U_1]\big(Z[U_1]B[U_1]-Z[U_2]B[U_2]\big),Z[U_1]\Lambda^s V\big)= B.$$
One has,
\begin{eqnarray*}
B=\big(\Lambda^s \big(Z[U_1]B[U_1]-Z[U_2]B[U_2]\big),\Lambda^s V\big)+\big(\big[\Lambda^s ,Z^{-1}[U_1]\big]Z[U_1]B[U_1]-Z[U_2]B[U_2],Z[U_1]\Lambda^s V\big)
\end{eqnarray*}
$B=B_1+B_2$.\\
\\
Now we have to estimate the terms ($B_1$) and ($B_2$).\\
\begin{eqnarray*}
(B_1)&=&\Big(\Lambda^s\big(\dfrac{-Q(\epsilon\zeta_1,\beta b,v_1)\beta \partial_x b g(\epsilon\zeta_1)v_1}{f(\epsilon\zeta_1)}+\dfrac{Q(\epsilon\zeta_2,\beta b,v_2)\beta \partial_x b g(\epsilon\zeta_2)v_2}{f(\epsilon\zeta_2)}\big),\Lambda^s \zeta_v\Big)\\
&\qquad+&\Big(\Lambda^s\big(\dfrac{\epsilon \gamma \beta q_1(\epsilon\zeta_1,\beta b)h_1(h_1+h_2)v_1^2\partial_x b}{(h_1+\gamma h_2)^3}-\dfrac{\epsilon \gamma \beta q_1(\epsilon\zeta_2,\beta b)h_1(h_1+h_2)v_2^2\partial_x b}{(h_1+\gamma h_2)^3}\big),\Lambda^s v \Big)
\end{eqnarray*}
With $Q(\epsilon\zeta_i,\beta b,v_i)=Q_0(\epsilon\zeta_i,\beta b)+\epsilon^2 Q_1(\epsilon\zeta_i,\beta b,v_i)$ for $i=1,2$.\\
\\
In order to control ($B_1$) we use the following decompositions,
\begin{multline}
\bullet \quad  \Big(\dfrac{ -Q_0(\epsilon\zeta_1,\beta b)\beta\partial_x b g(\epsilon\zeta_1)v_1}{f(\epsilon\zeta_1)}+\dfrac{ Q_0(\epsilon\zeta_2,\beta b)\beta\partial_x b g(\epsilon\zeta_2)v_2}{f(\epsilon\zeta_2)}\Big)\nn
\\=\Big(\dfrac{-Q_0(\epsilon\zeta_1,\beta b) g(\epsilon\zeta_1)}{f(\epsilon\zeta_1)}
+\dfrac{Q_0(\epsilon\zeta_2,\beta b)g(\epsilon\zeta_2)}{f(\epsilon\zeta_2)}\Big)(\beta \partial_x b v_1 )\nn\\
-\beta (v_1-v_2)\dfrac{Q_0(\epsilon\zeta_2,\beta b) g(\epsilon\zeta_2)\partial_x b}{f(\epsilon\zeta_2)}\nn.
\end{multline}

\begin{multline}
\bullet \quad  \beta\Big(\dfrac{-\epsilon^2 Q_1(\epsilon\zeta_1,\beta b,v_1)\partial_x b g(\epsilon\zeta_1)v_1}{f(\epsilon\zeta_1)}+\dfrac{\epsilon^2Q_1(\epsilon\zeta_2,\beta b,v_2)\partial_x b g(\epsilon\zeta_2)v_2}{f(\epsilon\zeta_2)}\Big)\nn
\\=\Big(\dfrac{-\epsilon^2 Q_1(\epsilon\zeta_1,\beta b,v_1) g(\epsilon\zeta_1)}{f(\epsilon\zeta_1)}
+\dfrac{\epsilon^2Q_1(\epsilon\zeta_2,\beta b,v_2) g(\epsilon\zeta_2)}{f(\epsilon\zeta_2)}\Big)(\beta\partial_x b v_1)\nn\\
-\beta(v_1-v_2)\dfrac{\epsilon^2Q_1(\epsilon\zeta_2,\beta b,v_2)g(\epsilon\zeta_2)\partial_x b }{f(\epsilon\zeta_2)}\nn.
\end{multline}

\begin{multline}
\bullet \quad \Big(\dfrac{\epsilon \gamma \beta q_1(\epsilon\zeta_1,\beta b)h_1(h_1+h_2)v_1^2\partial_x b}{(h_1+\gamma h_2)^3}-\dfrac{\epsilon \gamma \beta q_1(\epsilon\zeta_2,\beta b)h_1(h_1+h_2)v_2^2\partial_x b}{(h_1+\gamma h_2)^3}\Big)\\=\Big(\dfrac{\epsilon v_1 \gamma  q_1(\epsilon\zeta_1,\beta b)h_1(h_1+h_2)}{(h_1+\gamma h_2)^3}-\dfrac{\epsilon v_1 \gamma q_1(\epsilon\zeta_2,\beta b)h_1(h_1+h_2)}{(h_1+\gamma h_2)^3}\Big)(\beta \partial_x b v_1)\\+\Big(\dfrac{\epsilon \gamma  q_1(\epsilon\zeta_2,\beta b)h_1(h_1+h_2)\partial_x b}{(h_1+\gamma h_2)^3}\Big)\beta(v_1^2-v_2^2)\nn.
\end{multline}

Using the fact that, $\epsilon^2Q_1(\epsilon\zeta_i,\beta b,v_i)=Q_1(\epsilon\zeta_i,\beta b,\epsilon v_i)$, one deduces,

\begin{eqnarray*}
|B_1|&\leq& C(\beta|v_1|_{H^s},|b|_{H^{s+2}}) \epsilon |\zeta_1-\zeta_2|_{H^s}|\zeta_v|_{H^s}+ C(\epsilon|\zeta_2|_{H^s},|b|_{H^{s+2}})\beta|v_1-v_2|_{H^s}|\zeta_v|_{H^s}\\ & \quad +& C(\beta|v_1|_{H^s},|b|_{H^{s+1}})\epsilon |\zeta_1-\zeta_2|_{H^s}|\zeta_v|_{H^s}+ C(\epsilon|\zeta_2|_{H^s},\epsilon|v_2|_{H^s},|b|_{H^{s+1}})\beta|v_1-v_2|_{H^s}|\zeta_v|_{H^s}\\ & \quad +& C (\beta |v_1|_{H^s},\epsilon|v_1|_{H^s},|b|_{H^{s+1}})\epsilon|\zeta_1-\zeta_2|_{H^s}|v|_{H^s}\\ & \quad +& C (\epsilon|\zeta_2|_{H^s},\epsilon|v_1|_{H^s},\epsilon|v_2|_{H^s},|b|_{H^{s+1}})\beta |v_1-v_2|_{H^s}|v|_{H^s}.\\
&\leq&\max(\epsilon,\beta)  C_0 E^s(U_1-U_2)E^s(V).\\
&\leq&\max(\epsilon,\beta)  C_0 E^s(V)^2.
\end{eqnarray*}
with $C_0=C(M_{\rm CH},h^{-1},h_{03}^{-1},\max(\epsilon,\beta) \big\vert U_1\big\vert_{X^s},\max(\epsilon,\beta) \big\vert U_2\big\vert_{X^s},|b|_{H^{s+2}})$.\\
\\
The contribution of ($B_2$) is immediately bounded using Lemma~\ref{lem-nsZ}:
\begin{eqnarray*}
|B_2|&=&\big(\ \big[\Lambda^s,Z^{-1}[U_1]\big]\big(Z[U_1]B[U_1]-Z[U_2]B[U_2] \big)\ ,\ Z[U_1] \Lambda^s V\ \big)\\
&\leq& C|Z[U_1]B[U_1]-Z[U_2]B[U_2]|_{H^{s-1}\times H^{s-1}}|V|_{X^s}\\
&\leq& C\Big(\Big | \dfrac{-Q(\epsilon\zeta_1,\beta b,v_1)\beta \partial_x b g(\epsilon\zeta_1)v_1}{f(\epsilon\zeta_1)}+\dfrac{Q(\epsilon\zeta_2,\beta b,v_2)\beta \partial_x b g(\epsilon\zeta_2)v_2}{f(\epsilon\zeta_2)}\Big |_{H^{s-1}}\\
&\qquad +&\Big| \dfrac{\epsilon\gamma\beta q_1(\epsilon\zeta_1,\beta b)h_1(h_1+h_2)v_1^2\partial_x b}{(h_1+\gamma h_2)^3}-\dfrac{\epsilon\gamma\beta q_1(\epsilon\zeta_2,\beta b)h_1(h_1+h_2)v_2^2\partial_x b}{(h_1+\gamma h_2)^3}\Big |_{H^{s-1}}\Big)|V|_{X^s}\\
&\leq&\max(\epsilon,\beta)  C_0 E^s(U_1-U_2)E^s(V).\\
&\leq&\max(\epsilon,\beta)  C_0 E^s(V)^2.
\end{eqnarray*}
So we have,
$$|B|\leq C_0\max(\epsilon,\beta)E^s(V)^2.$$
Now one needs to control accordingly the right hand side $F$.\\
\\
In order to do so, we take advantage of the following Lemma.
\begin{lemma}\label{lemstZ}
Let $(\mu, \epsilon, \delta, \gamma, \beta, \bo)\in \P_{\rm CH}$ and $s\geq s_0>1/2$. Let $V=(\zeta_v,v)^\top$, $W=(\zeta_w,w)^\top\in X^{s}$ and $U_1=(\zeta_1,v_1)^\top$, $U_2=(\zeta_2,v_2)^\top\in X^{s}$, b $\in H^{s+2}$ such that there exists $h>0$ with
\[ 1-\epsilon\zeta_1\geq h>0,\ \quad 1-\epsilon\zeta_2\geq h>0,\ \quad \frac1\delta+\epsilon\zeta_1-\beta b\geq h>0,\ \quad \frac1\delta+\epsilon\zeta_2-\beta b\geq h>0.\]
Then one has
\begin{eqnarray*}\label{eq:est-dif-Z}
\Big\vert \ \Big(\ \Lambda^s \big( Z[U_1]-Z[U_2]\big) V \ , \ W \ \Big) \ \Big\vert \ & \leq & \ \epsilon\ C \  \big\vert U_1-U_2\big\vert_{X^s}\big\vert V\big\vert_{X^s}\big\vert W\big\vert_{X^0} \\
\label{eq:est-dif-ZA}
\Big(\ \Lambda^s \big( Z[U_1]A[U_1]-Z[U_2]A[U_2]\big) V \ , \ W \ \Big)  \ & \leq & \ \epsilon\ C\ \big\vert U_1-U_2\big\vert_{X^s}\big\vert V\big\vert_{X^s}\big\vert W\big\vert_{X^0}
\end{eqnarray*}
with $C=C(M_{\rm CH},h^{-1},\epsilon\big\vert U_1\big\vert_{X^s},\epsilon\big\vert U_2\big\vert_{X^s},|b|_{H^{s+2}})$.
\end{lemma}
\emph{Proof}.\\
We prove the Lemma~5.40 using the same techniques as in the Proof of~\cite[Lemma 7.2]{DucheneIsrawiTalhouk13}, adapted to our pseudo-symmetrizer as $\epsilon^2 Q_1(\epsilon\zeta_i,
\beta b,v_i)= Q_1(\epsilon\zeta_i,
\beta b,\epsilon v_i)$.\demo\\
\\
Let us continue the proof of Proposition~\ref{prop:stabilityZ}, by estimating $F$ defined in \eqref{rhsstabilityZ}.\\
\\
More precisely we want to estimate
\begin{eqnarray*}
\big(\ \Lambda^s F\ ,\ Z[U_1]\Lambda^s V\ \big) &= &\big(\ \Lambda^s F_1-\Lambda^s F_2\ ,\ Z[U_1]\Lambda^s V\ \big) \\
&\quad-&\big(\ \Lambda^s (\Sigma[U_1]-\Sigma[U_2])\partial_x U_2\ ,\ \Lambda^s V\ \big) \\
&\quad-&\big(\ \big[\Lambda^s,Z^{-1}[U_1]\big](\Sigma[U_1]-\Sigma[U_2])
\partial_xU_2\ ,\ Z[U_1]\Lambda^s V\ \big)\\
&\quad-&\big(\ \Lambda^s\big(Z[U_1]-Z[U_2]\big)(\partial_tU_2-F_2)\ ,\ \Lambda^s V\ \big)\\
&\quad-&\big(\ \big[\Lambda^s,Z^{-1}[U_1]\big](Z[U_1]-Z[U_2])(\partial_tU_2-F_2)\ ,\ Z[U_1]\Lambda^s V\ \big)\\&= &(I)+(II)+(III)+(IV)+(V).
\end{eqnarray*}
The contribution of $(I)$ is immediately bounded using Lemma~5.13. The contributions of $(II)$ and $(IV)$ follow Lemma~5.40. Finally, we control $(III)$ and $(V)$ using Lemma~5.13~\eqref{eq:est-com-Z-Z}.
All together, we proved using Lemma~5.12 that $F$ as defined in~\eqref{rhsstabilityZ}, satisfies
\begin{eqnarray*}\label{estFZ}
\vert (\Lambda^s F,Z[U_1]\Lambda^sV)\vert&\leq& \epsilon C\times( \vert \partial_x U_2\vert_{X^{s}}+ \vert \partial_t U_2-F_2 \vert_{X^{s}})E^s(V)^2+C E^s(V)E^s(F_1-F_2)\\&\leq& \max(\epsilon,\beta) C\times( \vert \partial_x U_2\vert_{X^{s}}+ \vert \partial_t U_2-F_2 \vert_{X^{s}})E^s(V)^2+C E^s(V)E^s(F_1-F_2).
\end{eqnarray*}
with $C=C(M_{\rm CH},h^{-1},h_{03}^{-1},\epsilon\big\vert U_1\big\vert_{X^s},\epsilon\big\vert U_2\big\vert_{X^s}, |b|_{H^{s+2}})$.\\
\\
Then one has
\begin{eqnarray*}\label{estFZ}
\vert (\Lambda^s F,Z[U_1]\Lambda^sV)\vert&\leq& \max(\epsilon,\beta) C_0\times( \vert \partial_x U_2\vert_{X^{s}}+ \vert \partial_t U_2-F_2 \vert_{X^{s}})E^s(V)^2+C_0 E^s(V)E^s(F_1-F_2).
\end{eqnarray*}
We can now conclude by Lemma~\ref{Lem:HsZ}, and the proof of Proposition~\ref{prop:stabilityZ}  is complete.\demo

\section{Full justification of the asymptotic model}\label{FJ}
A model is said to be {\em fully justified} (using the terminology of~\cite{Lannes}) if the Cauchy problem for both the full Euler system and the asymptotic model is well-posed for a given class of initial data, and over the relevant time scale; and if the solutions with corresponding initial data remain close.
We conclude our work by stating all the ingredients for the full justification of our model. Existence and uniqueness of the solution of the asymptotic model is given by
Theorem~6.1, while the stability with respect to the initial data is provided by Theorem~6.6.
\begin{theo}[Existence and uniqueness]\label{thbi1Z}
Let $\p=(\mu, \epsilon, \delta, \gamma, \beta, \bo)\in \P_{\rm CH}$ and  $s\geq s_0+1$, $s_0>1/2$, and assume
	 $U_0=(\zeta_0,v_0)^\top\in X^s$, $b\in H^{s+2}$  satisfies~\eqref{CondDepth},\eqref{CondEllipticity}, and \eqref{H3}.
Then there exists a maximal time $T_{\max}>0$, uniformly bounded from below with respect to $\p\in \P_{\rm CH}$, such that the system of
equations~\eqref{eq:Serre2} admits a unique strong solution $U=(\zeta,v)^\top \in C^0([0,T_{\max});X^s)\cap C^1([0,T_{\max});X^{s-1})$ with the initial value $(\zeta,v)\id{t=0}=(\zeta_0,v_0)$, and preserving the conditions~\eqref{CondDepth},\eqref{CondEllipticity} and \eqref{H3} (with different lower bounds) for any $t\in [0,T_{\max})$.\\
\\
Moreover, there exists $\lambda,C_0 = C(h_{01}^{-1},h_{02}^{-1},h_{03}^{-1},M_{\rm CH},T,\big\vert U_0\big\vert_{X^{s}},|b|_{H^{s+2}})$, independent of $\p\in \P_{\rm CH}$, such that $T_{\max}\geq T/\max(\epsilon,\beta)$, and one has the energy estimates\\
   \\
$\forall\ 0\leq t\leq\frac{T}{\max(\epsilon,\beta)}$,
\[\big\vert U(t,\cdot)\big\vert_{X^{s}} \ + \
\big\vert \partial_t U(t,\cdot)\big\vert_{X^{s-1}}  \leq C_0e^{\max(\epsilon,\beta)\lambda t}+\max(\epsilon,\beta) C_0^2 \int_0^t e^{\max(\epsilon,\beta)\lambda(t-t')}dt'+\max(\epsilon,\beta) C_0
\]
If $T_{\max}<\infty$, one has
         \[ \vert U(t,\cdot)\vert_{X^{s}}\longrightarrow\infty\quad\hbox{as}\quad t\longrightarrow T_{\max},\]
         or one of the conditions~\eqref{CondDepth},\eqref{CondEllipticity}, \eqref{H3} ceases to be true as $ t\longrightarrow T_{\max}$.

\end{theo}
\emph{Proof}.\\
We construct a sequence of approximate solution $(U^n=(\zeta^n,v^n))_{n\ge 0}$ through the
induction relation
\begin{equation}\label{approximatesysz}
         U^0=U_0,\quad\mbox{ and }\quad
	\forall n\in\NN, \quad
	\left\lbrace
	\begin{array}{l}
	\dsp\partial_t U^{n+1}+A[U^n]\partial_x U^{n+1}+B[U^n]
	=0;
        \\
	\dsp U^{n+1}_{\vert_{t=0}}=U_0.
	\end{array}\right.
\end{equation}
\\
By Proposition~\ref{ESpropZ}, there exists $U^{n+1}\in C^0([0,\dfrac{T_{n+1}}{\max(\epsilon,\beta)}];X^s)\cap C^1([0,\dfrac{T_{n+1}}{\max(\epsilon,\beta)}];X^{s-1})$ unique solution to~\eqref{approximatesysz} if $U^{n}\in C^0([0,\dfrac{T_n}{\max(\epsilon,\beta)}];X^s)\cap C^1([0,\dfrac{T_n}{\max(\epsilon,\beta)}];X^{s-1})\subset X^s_T$, and satisfies~\eqref{CondDepth},\eqref{CondEllipticity} and \eqref{H3}.\\
\\
\noindent {\em Existence and uniform control of the sequence $U^n.$}\\
\\
The existence of $T'>0$ such that the sequence $U^n$ is  uniquely defined, controlled in $X^s_{T'}$, and satisfies~\eqref{CondDepth},\eqref{CondEllipticity} and \eqref{H3}, uniformly with respect to $n\in\NN$, is obtained by induction, as follows.

Proposition~\ref{ESpropZ} yields
\begin{eqnarray*}\label{energyestimates-n} E^s(U^{n+1})(t) \leq  e^{\max(\epsilon,\beta) \lambda_n t}E^s(U_0)+\max(\epsilon,\beta)C_n\int_0^t e^{\max(\epsilon,\beta)\lambda_n (t-t')}dt'.\end{eqnarray*}
\begin{eqnarray*}\big\vert \partial_t U^{n+1}(t,\cdot)\big\vert_{X^{s-1}} &\leq& C_n E^s(U^{n+1})(t)+\max(\epsilon,\beta)C_n\\
&\leq& C_n e^{\max(\epsilon,\beta) \lambda_n t}E^s(U_0)+\max(\epsilon,\beta)C_n^2\int_0^t e^{\max(\epsilon,\beta)\lambda_n (t-t')}dt'+\max(\epsilon,\beta)C_n ,\end{eqnarray*}
with  $C_n,\lambda_n=C(M_{\rm CH},h_{01,n}^{-1},h_{02,n}^{-1},T_n,\big\Vert U^n\big\Vert_{X^s_{T_n}},|b|_{H^{s+2}})$, provided $U^n\in X^s_{T_n} $ satisfies~\eqref{CondDepth},\eqref{CondEllipticity} and \eqref{H3}  with positive constants $h_{01,n},h_{02,n}$, and $h_{03,n}$ on $[0,T_n/\max(\epsilon,\beta)]$.\\
\\
It is a consequence of the work~\cite[Theorem 7.3]{DucheneIsrawiTalhouk13} and by taking into account the topographic variation that the assumptions \eqref{CondDepth} and \eqref{CondEllipticity} may be imposed only on the initial data and then is automatically satisfied over the relevant time scale.
Let us prove it now for \eqref{H3},\\
\\
Since $U^{n}=(\zeta^n,v^n)^\top$ satisfies~\eqref{approximatesysz}, one has
\begin{eqnarray*}
\partial_t \zeta ^{n+1}=-f(\epsilon\zeta^n)\partial_x v^{n+1}-\epsilon f'(\epsilon\zeta^n)v^n\partial_x \zeta^{n+1}+\beta\partial_x b g(\epsilon \zeta ^n) v^n,
\end{eqnarray*}
and
\begin{eqnarray*}
\partial_t v^{n+1}&=& -\mathfrak{T}[\epsilon\zeta^n,\beta b]^{-1}\Big(Q_0(\epsilon\zeta^n,\beta b)\partial_x\zeta^{n+1}+\epsilon \mathfrak{Q}[\epsilon\zeta^n,\beta b,v^n]\partial_x v^{n+1}+\epsilon^2Q_1(\epsilon\zeta^n,\beta b,v^n)\partial_x\zeta^{n+1} \\ &+& \epsilon \dfrac{\gamma\beta q_1(\epsilon\zeta^n,\beta b)h_1(h_1+h_2)v^{n^2}\partial_x b}{(h_1+\gamma h_2)^3} \Big)-\epsilon\varsigma v^n\partial_x v^{n+1}.
\end{eqnarray*}

Using continuous Sobolev embedding of $H^{s-1}$ into $L^\infty$ ($s-1>1/2$), and since $U^n$ satisfies~\eqref{CondDepth},\eqref{CondEllipticity} with $h_{01,n},h_{02,n}$ on $[0,T_n/\max(\epsilon,\beta)]$, one deduces that
 \begin{equation}\label{zetat}
\vert \partial_t \zeta^{n+1} \vert_{L^{\infty}}
\leq C(M_{\rm CH},h_{01,n}^{-1},h_{02,n}^{-1},\beta|b|_{H^s}\big)\big\Vert U^n\big\Vert_{X^s_{T_n}},
\end{equation}
and
 \begin{equation}\label{vt}
\vert \partial_t v^{n+1} \vert_{L^{\infty}}
\leq C(M_{\rm CH},h_{01,n}^{-1},h_{02,n}^{-1},\beta|b|_{H^{s+1}}\big)\big\Vert U^n\big\Vert_{X^s_{T_n}}.
\end{equation}
Let $g^{n+1}=a_1(\epsilon \zeta^{n+1},\beta b)+a_2(\epsilon \zeta^{n+1},\beta b)\epsilon^2(v^{n+1})^2$,\\ where $\big(a_1(\epsilon \zeta^{n+1},\beta b),a_2(\epsilon \zeta^{n+1},\beta b)\big)=\Big((\gamma+\delta)q_1(\epsilon \zeta^{n+1},\beta b)-\mu\beta \omega \partial_x^2 b,-\gamma q_1(\epsilon \zeta^{n+1},\beta b)\dfrac{(h_1+h_2)^2}{(h_1+\gamma h_2)^3}\Big)$ \\
One has,
\begin{eqnarray*}
g^{n+1}&=&g^{n+1}\id{t=0}+\int_0^t \partial_t a_1(\epsilon \zeta^{n+1},\beta b)+\epsilon^2\int_0^t \partial_t a_2(\epsilon \zeta^{n+1},\beta b)(v^{n+1})^2\\&\qquad+& 2\epsilon^2\int_0^t a_2(\epsilon \zeta^{n+1},\beta b)v^{n+1}\partial_tv^{n+1}\\
&=&g^{n+1}\id{t=0}+(\gamma+\delta)\epsilon\kappa_1\int_0^t \partial_t\zeta^{n+1}+\epsilon^3\int_0^ta'_2(\epsilon\zeta^{n+1},\beta b)\partial_t\zeta^{n+1} (v^{n+1})^2\\&\qquad+&2\epsilon^2\int_0^t a_2(\epsilon \zeta^{n+1},\beta b)v^{n+1}\partial_tv^{n+1}
\end{eqnarray*}
so that~\eqref{zetat}  and~\eqref{vt} yields
\begin{eqnarray*}\big\vert g^{n+1}-g^{n+1}\id{t=0}\big\vert_{L^\infty}\leq \epsilon t \times C(M_{\rm CH},h_{01,n}^{-1},h_{02,n}^{-1},|b|_{H^{s+2}}\big)\big\Vert U^n\big\Vert_{X^s_{T_n}}.\end{eqnarray*}
Now, one has $g^{n+1}\id{t=0}\equiv g^0\id{t=0} \ge h_{03,0}>0$, independent of $n$. Thus one can easily prove (by induction) that it is possible to chose $T'>0$ such that $g^{n+1}>\dfrac{h_{03}}{2}$ holds on $[0,T'/\max(\epsilon,\beta)]$, and the above energy estimates hold uniformly with respect to $n$, on $[0,T'/\max(\epsilon,\beta)]$.

More precisely, one has that
$U^n$ satisfies \eqref{H3} with $\dfrac{h_{03}}{2}>0$
  and the estimates
\begin{eqnarray*}\label{energyestimates-nz} E^s(U^{n})(t) \leq  e^{\max(\epsilon,\beta)\lambda t}E^s(U_0)+\max(\epsilon,\beta) C_0\int_0^t e^{\max(\epsilon,\beta)\lambda(t-t')}dt'\end{eqnarray*}
\begin{eqnarray}\big\vert \partial_t U^{n}(t,\cdot)\big\vert_{X^{s-1}} &\leq& C_0 E^s(U^{n})(t)+\max(\epsilon,\beta)C_0\nn\\&\leq& C_0 e^{\max(\epsilon,\beta)\lambda t}E^s(U_0)+\max(\epsilon,\beta) C_0 ^2\int_0^t e^{\max(\epsilon,\beta)\lambda(t-t')}dt'+\max(\epsilon,\beta)C_0.\nn\\\end{eqnarray}
on $[0,T'/\max(\epsilon,\beta)]$, where $C_0,\lambda=C(M_{\rm CH},h_{01}^{-1},h_{02}^{-1},h_{03}^{-1},T',\big\vert U_0\big\vert_{X^s},|b|_{H^{s+2}})$ are uniform with respect to $n$.\\
\\
For the completion of the proof ({\em Convergence of $U^n$ towards a solution of the nonlinear problem}) we use
the same techniques as in the proof of~\cite[Theorem 7.3]{DucheneIsrawiTalhouk13}(see e.g.~\cite{AlinhacGerard}). The uniqueness of $U$ follows from the priori estimate result of Proposition~\ref{prop:stabilityZ} with $F_1\equiv F_2\equiv 0$. \demo

\begin{theo}[Stability] \label{th:stabilityWP}
Let $\p=(\mu, \epsilon, \delta, \gamma, \beta, \bo)\in \P_{\rm CH}$ and $s\geq s_0+1$ with $s_0>1/2$, and assume
	 $U_{1,0}=(\zeta_{1,0},v_{1,0})^\top\in X^{s}$, $U_{2,0}=(\zeta_{2,0},v_{2,0})^\top\in X^{s+1}$, and $b\in H^{s+2}$ satisfies~\eqref{CondDepth},\eqref{CondEllipticity}, and \eqref{H3}. Denote $U_j$ the solution to~\eqref{eq:Serre2} with $U_j\id{t=0}=U_{j,0}$. \\Then there exists $T,\lambda,C_0= C(M_{\rm CH}, h_{01}^{-1},h_{02}^{-1},h_{03}^{-1},\big\vert U_{1,0}\big\vert_{X^s},\vert U_{2,0}\vert_{X^{s+1}},|b|_{H^{s+2}})$ such that $\forall t\in [0,\frac{T}{\max(\epsilon,\beta)}]$,
 \begin{equation*}
\big\vert (U_1-U_2)(t,\cdot)\big\vert_{X^s} \leq C_0 e^{\max(\epsilon,\beta)\lambda t} \big\vert U_{1,0}-U_{2,0}\big\vert_{X^s}.
\end{equation*}
\end{theo}
\emph{Proof}.\\
The existence and uniform control of the solution $U_1$ (resp. $U_2$) in $L^\infty([0,T/\max(\epsilon,\beta)];X^s)$ (resp. $L^\infty([0,T/\max(\epsilon,\beta)];X^{s+1})$) is provided by Theorem~\ref{thbi1Z}.
The proposition is then a direct consequence of the {\em a priori} estimate of Proposition~\ref{prop:stabilityZ}, with $F_1=F_2=0$, and Lemma~\ref{lemmaes}.\demo
\\

Finally, the following ``convergence result" states that the solutions of our system approach the solutions of the full Euler system, with as good a precision as $\mu$  is small.

\begin{theo}[Convergence]\label{th:convergence}
Let $\p=(\mu, \epsilon, \delta, \gamma, \beta, \bo)\in \P_{\rm CH}$(see~\eqref{eqn:defRegimeCHmr}) and $s\geq s_0+1$ with $s_0>1/2$, and let $U^0\equiv(\zeta^0,\psi^0)^\top\in H^{s+N}(\RR)^2$, $ b \in H^{s+N}$ satisfy the hypotheses~\eqref{CondDepth},\eqref{CondEllipticity}, and \eqref{H3}, with $N$ sufficiently large. We suppose  $U\equiv (\zeta,\psi)^\top$ a unique solution to the full Euler system~\eqref{eqn:EulerCompletAdim} with initial data $(\zeta^0,\psi^0)^\top$, defined on $[0,T_1]$ for $T_1>0$ \footnote[1]{ To our knowledge, the local well-posedness of the full Euler system in the two-fluid configuration over a variable topography seems to be an open problem.}, and we suppose that $U\equiv(\zeta,\psi)^\top$ satisfies the assumptions of our consistency result, Proposition~\ref{th:ConsSerreVar}. Then there exists $C,T>0$, independent of $\p$, such that

\begin{itemize}
\item There exists a unique solution $U_a\equiv (\zeta_a,v_a)^\top$ to our new model~\eqref{eq:Serre2}, defined on $[0,T]$ and with initial data $(\zeta^0,v^0)^\top$ (provided by Theorem~\ref{thbi1Z});
\item With $v$, defined as in~\eqref{defv}, one has for any $t\in[0,T]$,
\[ \big\vert (\zeta,v)-(\zeta_a,v_a) \big\vert_{L^\infty([0,t];X^s)}\leq C \ \mu^2\ t.\]
\end{itemize}
\end{theo}
\emph{Proof}.\\
The existence of $U_a$ is given by our Theorem~\ref{thbi1Z} (we choose $T$ as the minimum of the existence time of both solutions; it is bounded from below, independently of $\p\in\P_{\rm CH}$). Assuming that $U\equiv(\zeta,\psi)^\top$ satisfies the assumptions of our consistency result, Proposition~\ref{th:ConsSerreVar}, therefore $(\zeta, v)^\top$  solves~\eqref{eq:Serre2} up to a residual $R=(r_1,r_2)^\top$, with $\big\vert R\big\vert_{L^\infty([0,T];H^s)}\leq C(M_{\rm SW},h_{01}^{-1},|b|_{H^{s+N}},\big\vert U^0\big\vert_{H^{s+N}})(\mu^2+\mu\epsilon^2+\mu\beta^2+\mu\epsilon\beta)$. As a matter of fact, since $\p\in \P_{\rm CH}$ therefore the residual is now bounded by $\mu^2$.  The result follows from the stability Proposition~\ref{prop:stabilityZ}, with $F_1=(r_1,\mfT[\epsilon\zeta]^{-1}r_2)^\top$ and $F_2=0$.\demo

\end{document}